\documentclass{amsart}
\usepackage{fullpage}
\usepackage[utf8]{inputenc}
\usepackage[T1]{fontenc}
\usepackage{lmodern}
\usepackage{amsmath,amsfonts,amssymb,amsthm,mathtools,mathrsfs}
\usepackage{xcolor}
\usepackage{enumitem}
\usepackage[hidelinks]{hyperref}
\allowdisplaybreaks
\hypersetup{pdftitle={Strong convergence rates for numerical approximations of the stochastic Burgers equation with additive rough noise},pdfauthor={Ziyi Lei and Houqi Su}}
\newtheorem{theo}{Theorem}[section]
\newtheorem{lemma}[theo]{Lemma}

\theoremstyle{remark}
\newtheorem{rem}[theo]{Remark}
\newcommand{\N}{\mathbb N}
\newcommand{\R}{\mathbb R}
\newcommand{\E}{\mathbb E}
\newcommand{\HH}{\mathcal H}
\newcommand{\IF}{\mathcal F}
\newcommand{\IL}{\mathcal L}
\newcommand{\ID}{\mathcal D}

\makeatletter
\newcommand{\printaddressesonfirstpage}{%
 \ifx\@empty\addresses\else
 \insert\footins{\footnotesize\interlinepenalty\@M\@setaddresses}%
 \global\let\addresses\@empty\fi}
\makeatother

\title{Strong convergence rates for numerical approximations of the stochastic Burgers equation with additive rough noise}
\author{Ziyi Lei}
\address{Department of Mathematics, College of Sciences, China University of Mining and Technology, Beijing 100083, China}
\email{ziyilei@lsec.cc.ac.cn}
\author{Houqi Su$^*$}
\thanks{$^*$ Corresponding author.}
\address{School of Mathematical Sciences, Capital Normal University, Beijing 100048, China}
\email{marksu@amss.ac.cn}

\begin{document}

\begin{abstract}
We establish strong convergence rates for a fully discrete approximation of the one-dimensional stochastic Burgers equation driven by additive rough noise \((-A)^\alpha dW\), \(0\le\alpha<1/4\), below the critical roughness threshold $\alpha=1/4$. Through the Cole--Hopf transformation, the nonlinear equation is converted into a linear stochastic heat equation with spatially correlated multiplicative noise, which is discretized by a spectral Galerkin method in space and an accelerated exponential Euler method in time. A numerical approximation of the original solution is then reconstructed through a scale-invariant regularized inverse Cole--Hopf transformation. The convergence analysis combines strong Sobolev error estimates with strict positivity and uniform inverse-moment bounds for the exact transformed solution. Under an initial approximation condition and a weak-\(L^4\) tail assumption on the transformed initial datum, we obtain the spatial strong convergence rates arbitrarily close to \(1/2-2\alpha\) and a temporal rate of order \(1/2\), up to a logarithmic correction. For $0\le\alpha<1/4$, the spatial rate is essentially sharp within the spectral 
Galerkin framework, the temporal rate is independent of the roughness parameter 
$\alpha$, and we prove a converse estimate showing both the necessity of the 
initial approximation condition and the weak-$L^4$ tail of the initial datum.

\end{abstract}
\maketitle
\printaddressesonfirstpage

\section{Introduction}

The stochastic Burgers equation (SBE) is a canonical model for the interplay between nonlinear transport, viscous dissipation, and random forcing. Its deterministic counterpart is a prototypical nonlinear advection--diffusion equation, in which the nonlinear transport tends to generate steep gradients while viscosity provides a regularizing mechanism \cite{bec2007burgers}. The stochastic forcing incorporates random fluctuations and unresolved effects, leading to a rich combination of nonlinear and probabilistic phenomena \cite{bec2007burgers,bertini1997stochastic}. The SBE and its variants arise in a variety of applications, including fluid dynamics and turbulence \cite{bec2007burgers}, traffic flow \cite{nagatani1998burgers}, and, through their close relation to the Kardar--Parisi--Zhang (KPZ) equation and directed polymers, interface growth and disordered systems \cite{bertini1997stochastic,halpin1995kinetic}. In this paper, we study the one-dimensional SBE driven by additive rough noise. Let $H:=L^2(0,1;\mathbb R)$ and let \(A\) denote the Dirichlet Laplacian on $H$ with domain \(H^2(0,1)\cap H_0^1(0,1)\). For \(T>0\) and \(0\le\alpha<1/4\), we consider 
\begin{equation}\label{eq:sharp-sbe}
\begin{cases}
du(t)=Au(t)\,dt+\partial_x\bigl(u(t)^2\bigr)\,dt+(-A)^{\alpha}dW(t),& t\in(0,T],\\
u(0)=u_0,
\end{cases}
\end{equation}
where \(u_0\) is an \(\mathcal F_0\)-measurable \(H\)-valued random variable, and \(W\) is a cylindrical Wiener process on \(H\). The precise probabilistic setting and assumptions are given in Subsection~\ref{sec:2.1}.

The analytical theory of the SBE, particularly its well-posedness, 
has been extensively investigated in the literature; see, for example, \cite{bertini1994stochastic,daprato1994stochastic,daprato1995stochastic,gyongy1998existence}. For the Dirichlet problem on \((0,1)\), the well-posedness under suitable assumptions was established in \cite{daprato1994stochastic}. On the whole real line, Bertini et al. \cite{bertini1994stochastic} studied the equation through the Cole--Hopf transformation and obtained a probabilistic representation of the solution via a generalized Feynman--Kac formula. Despite these analytical advances, 
numerical approximations remain important for studying the solution dynamics and computing statistical quantities. 
Existing convergence analyses for the SBE have primarily addressed the space--time white noise (\(\alpha=0\)), or spatially smoother additive noise. For the space--time white noise, spatial semi-discretizations based on finite difference methods (cf. \cite{hairer2011approximations,alabert2006numerical,hairer2012spatial,hairer2014approximating}) and spectral Galerkin methods (cf. \cite{blomker2013galerkin}) have been analyzed, with pathwise convergence rates obtained under suitable assumptions (see also \cite{hairer2016optimal} for optimal pathwise rates for Burgers-type equations). Subsequently, the temporal and fully discrete approximations have  received increasing attention. For instance, Khan \cite{khan2021higher} proposed higher-order pathwise approximations aimed at improving temporal convergence rates, while Khan et al. \cite{khan2024simplified} studied a simplified exponential Euler scheme and derived corresponding temporal error estimates. More recently, Jentzen et al. \cite{jentzen2019strong} proved strong convergence for a fully discrete, nonlinearity-truncated accelerated exponential Euler-type approximation, and Mazzonetto et al. \cite{mazzonetto2020existence} developed a unified framework covering both almost sure and strong convergence of numerical approximations. By contrast, 
when \(0<\alpha<1/4\), the driving noise of the SBE \eqref{eq:sharp-sbe}  is spatially rougher than the space–time white noise, which further reduces the spatial regularity of the solution. Consequently, the spatial regularity of the stochastic convolution is limited to Sobolev orders below \(1/2-2\alpha\), restricting the attainable spatial convergence rates. Moreover, the reduced regularity, together with the derivative structure and non-global Lipschitz continuity of the nonlinear drift, makes the stability and error estimates for the  numerical approximations more challenging. To the best of our knowledge, quantitative strong convergence rates for fully discrete approximations of the SBE \eqref{eq:sharp-sbe} in this rough-noise regime have not been established prior to the present work.

To fill this gap, we develop a fully discrete approximation that achieves a spatial strong convergence rate arbitrarily close to \(1/2-2\alpha\) and a temporal rate of order \(1/2\), up to a logarithmic correction. 
 Following the Cole–Hopf approach in \cite{bertini1994stochastic}, we formally transform the nonlinear SBE  \eqref{eq:sharp-sbe} into a stochastic heat equation of parabolic Anderson type with multiplicative spatially correlated noise: 
\begin{equation}\label{SPDE:HC-transform-SEE-I}
\begin{aligned}
U(t)
={}&S_1(t)U_0
+\int_0^t S_1(t-s)\rho_\alpha U(s)\,ds +\int_0^t S_1(t-s)G_{\alpha}(U(s))\,d\beta(s),
\qquad t\in(0,T]
\end{aligned}
\end{equation}
with $U(0):=\exp\bigl(\int_0^{\cdot}u_0(y)\,dy-\min_{z\in[0,1]}\int_0^zu_0(y)\,dy\bigr).$ The coefficients and the driving noise in \eqref{SPDE:HC-transform-SEE-I} are specified in Subsection \ref{sec:2.2}. Exploiting the linear structure of the transformed equation, we construct a fully discrete scheme by utilizing a spectral Galerkin method in space and an accelerated exponential Euler method in time:
\begin{align}\label{solution:HC-transform-SEE-fulldiscrete-I}
   U^N_{m}=S_1(\tau)U^N_{m-1}+ \int_{t_{m-1}}^{t_m}S_1(t_m-s)\Pi_N(\rho_\alpha U^N_{m-1})\,ds+\sum_{k\in\N_+}\!\!\int_{t_{m-1}}^{t_m}S_1(t_m-s)\Pi_N(U^N_{m-1}b_{k,\alpha})\,d\beta_k(s)
\end{align}
for all $m,N\in\mathbb{N}_+$ with initial value $U^N_0:=\Pi_N U_0$. Although the exact transformed solution \(U\) is almost surely strictly positive on \([0,T]\times[0,1]\), fully discrete approximation \(U_m^N\) does not preserve the strict positivity. To recover an approximation of the original Burgers solution, we introduce a scale-invariant relative regularization of the logarithmic derivative. Set $\delta:=1\wedge(N^{-\varsigma}+\sqrt{\tau})$ for a fixed $\varsigma\in(0,1/2-2\alpha)$. Let \(P_N\) be the orthogonal projection onto the first \(N\) eigenfunctions of the Dirichlet Laplacian. We introduce the mapping
\[
\mathscr R_{N,\delta}(f)
:=
\begin{cases}
\displaystyle
P_N\Big(
\frac{\partial_x f}
{f\vee\bigl(\delta\|f\|_{\mathcal H_1^1}\bigr)}
\Big),
& f\neq0,\\[3mm]
0,
& f=0.
\end{cases}
\]
This relative regularization yields a uniform bound in \(N\) and the input function on the inverse transformation and, together with the strict positivity and inverse-moment estimates of $U$, enables the approximation error transfer to the original Burgers solution. The fully discrete approximation to the original SBE~\eqref{eq:sharp-sbe} is defined by\begin{align}\label{fulldiscrete-u-I}
u_m^{N,\delta}
:=
\mathscr R_{N,\delta}(U_m^N),
\qquad m=1,\ldots,M,
\end{align}
with the initial approximation $u_0^{N,\delta}
:=
\mathcal R_{N,N^{-\varsigma}}(U_0^N)$. Under the assumptions of Theorem \ref{thm:sharp-main}, the resulting approximation satisfies the following strong error estimate, with a constant
\(C_{\alpha,\varsigma,T}>0\) independent of \(N\) and \(\tau\): 
\begin{align}
     \sup_{0\le m\le M}\big\|u(t_m)-u^{N,\delta}_m\big\|_{L^2(\Omega;H)}\le C_{\alpha,\varsigma,T}\big(1+\mathcal{J}_{\varsigma}+\big\|\|U_0\|_{\mathcal{H}_1^1}\big\|_{4,\infty}^4\big)^{1/2}\big(N^{-\varsigma}+\tau^{1/2}(1+0\vee\ln (1/\tau))^{1/2}\big),
 \end{align}   
where $\mathcal{J}_{\varsigma}:=\sup_{N\ge1}N^{2\varsigma}
 \E\big[\big\|u_0-
 \mathcal R_{N,N^{-\varsigma}}(\Pi_NU_0)\big\|_H^2\big]$ and $ \|\cdot\|_{4,\infty}^4:=\sup_{R\ge1}R^4\mathbb P(|\cdot|>R)$.

The main contributions of this work, together with the key analytical challenges addressed, are summarized as follows:

\begin{enumerate}[
    label=(\arabic*),
    leftmargin=2.0em,
    labelsep=0.6em,
    itemsep=0.4\baselineskip,
    topsep=0.4\baselineskip,
    parsep=0.4\baselineskip
]

\item \textbf{Sharp spatial and
$\boldsymbol{\alpha}$-independent temporal convergence rates.}

The spatial convergence rate
$\varsigma<1/2-2\alpha$ reflects the limited spatial
regularity induced by the high-frequency components of the
noise. For $0<\alpha<1/4$, it is essentially sharp within
the spectral Galerkin framework
(see Remark~\ref{rem:sharp-Galerkin}). The main challenge in the temporal analysis arises from the limited temporal regularity of the original Burgers solution, which makes direct regularity-based error estimates insufficient to obtain the desired half-order convergence rate. By exploiting the Cole--Hopf transformation and an accelerated exponential Euler discretization, together with the Fourier multiplier estimates (Lemma~\ref{lem:sharp-multiplier}), we control the temporal discretization error through weighted estimates of solution increments in weaker fractional Sobolev spaces \(\mathcal H_1^\gamma\), where \(0<\gamma<1/2-2\alpha\). These estimates also account for the initial-layer singularity and lead to a temporal strong convergence rate \(1/2\), up to a logarithmic correction, with a convergence exponent independent of \(\alpha\). 


\item \textbf{Strict positivity and stable inverse Cole--Hopf reconstruction.}

A major difficulty in transferring strong error estimates back to the original Burgers solution is controlling the division by a random denominator, which may become arbitrarily small or non-positive after numerical discretization. At the continuous level, we establish the strict positivity and uniform inverse-moment bounds for the exact transformed solution $U$ through a finite-noise reflected Brownian motion Feynman--Kac formula. At the discrete level, we design a relative regularization that truncates the denominator. This preserves the natural scale invariance $U \mapsto cU$ of the inverse Cole--Hopf mapping and furnishes a deterministic stability bound $\Vert{}\mathscr{R}_{N,\delta}(f)\Vert{}_H \le \delta^{-1}$. These ingredients, combined with quantitative error propagation estimates, enable stable reconstruction of the original Burgers solution from the numerical approximation of the transformed equation.

\item \textbf{Weak initial-data assumptions and necessary tail bounds.}

To handle both the low regularity of $u_0 \in H$ and the heavy tails of $U_0$ without requiring finite fourth moments, we combine an initial-layer splitting—pairing parabolic smoothing outside the layer with deterministic Burgers comparisons (Lemma \ref{lem:deterministic-HC}) inside—with $\mathcal{F}_0$-conditional localization and layer-cake tail estimates. This yields uniform strong convergence under two minimal conditions: $\mathcal{J}_\varsigma < \infty$ and $\|\|U_0\|_{\mathcal{H}_1^1}\|_{4,\infty} < \infty$. Crucially, we prove a converse theorem (Remark \ref{rem:sharp-Galerkin}(2)) demonstrating that any initial law achieving the asserted convergence rate must satisfy $\mathcal{J}_\varsigma < \infty$ and $\|\|u_0\|_H\|_{4,\infty} < \infty$. Thus, $\mathcal{J}_\varsigma < \infty$ is strictly necessary and the weak-$L^4$ tail exponent $4$ is sharp for $L^2(\Omega;H)$ convergence, with the sufficient and necessary conditions coinciding whenever $\int_0^\cdot u_0(y)\,dy$ is uniformly bounded.


\end{enumerate}

For a broader context, related questions have been studied under various discretization and noise settings. In the space-time white-noise setting, Alabert and Gy\"ongy \cite{alabert2006numerical} analyzed finite difference approximations and obtained a spatial convergence rate arbitrarily close to \(1/2\) in the \(L^2(0,1)\)-norm. For spatially smoother noise, a broader range of convergence results is available. Bl\"omker et al. \cite{blomker2013full} studied a fully discrete approximation of the SBE driven by spatially correlated \(Q\)-Wiener noise and derived pathwise error estimates. Under additive trace-class noise, Hutzenthaler et al. \cite{hutzenthaler2019rates} established strong convergence rates for fully discrete approximations, while Hutzenthaler and Link \cite{hutzenthaler2022multiplicative} obtained corresponding results for the multiplicative trace-class noise. Weak approximation has also been studied in this setting: Br\'ehier et al. \cite{brehier2026weak} derived weak error rates for spectral Galerkin approximations. Related convergence results are also available for other Burgers-type equations. In particular, Wang et al. \cite{wang2026burgershuxley} established strong spatial and temporal convergence rates for the stochastic Burgers--Huxley equation driven by additive space-time white noise.

The remainder of the paper is organized as follows. Section~\ref{sec:method} introduces the preliminary setting, formulates the Cole--Hopf transformation, constructs the fully discrete scheme for \eqref{eq:sharp-sbe}, and states the main convergence result. Section~\ref{sec:linear} establishes the well-posedness, strict positivity, inverse-moment bounds, and strong \(\mathcal H_1^\gamma\)-error estimates for the transformed linear equation. Section~\ref{sec:reconstruction} proves the nonlinear reconstruction stability lemma and combines the weak-tail localization with the initial-layer analysis to complete the proof of Theorem~\ref{thm:sharp-main}. Section~5 presents numerical experiments that corroborate the theoretical convergence rates. Finally, auxiliary analytic, semigroup, and Cole--Hopf identities are collected in Appendices~\ref{app:analytic} and~\ref{app:colehopf}.

\section{The set-up, numerical method, and main result}\label{sec:method}


\subsection{Notation and assumptions}\label{sec:2.1}

Throughout this paper, $H$ denotes the Hilbert space $L^2(0, 1; \mathbb{R})$ equipped with the usual norm $\Vert{}\cdot\Vert{}_H$. More generally, for any $p \in [1, \infty]$, we use $\Vert{}\cdot\Vert{}_p$ to denote the norm in $L^p(0, 1; \mathbb{R})$. We denote by $\ell^2$ the Hilbert space of all square-summable real sequences.

Let $\IL(U_1;U_2)$ and $\IL_2(U_1;U_2)$ be the spaces of bounded linear operators and Hilbert--Schmidt operators, respectively, from separable Hilbert spaces $U_1$ to $U_2$. 
In particular, we use the notations $\IL(U_1):=\IL(U_1;U_1)$ and $\IL_2(U_1):=\IL_2(U_1;U_1)$. Additionally, we denote by $C([0,1]; \mathbb{R})$ the space of continuous real-valued functions on $[0,1]$, equipped with the supremum norm  $\Vert{}\cdot\Vert{}_{C([0,1];\R)}$. For  $a \in \mathbb{N}$ and $p \in [1, \infty]$, we denote the standard Sobolev space on $(0, 1)$ by $W^{a,p}(0, 1)$ with \(W^{0,p}(0,1):=L^p(0,1;\R)\); we refer to \cite[Section 1.11.2]{yagi2010abstract} for its precise definition.
In particular, the zero-trace space
used below is defined by
$
 W^{1,2}_0:=\{f\in W^{1,2}(0,1)\colon f(0)=f(1)=0\},
$
where the endpoint values are understood in the sense of the continuous representative. The dual space $(W_0^{1,2}(0, 1))^*$ consists of all bounded linear functionals on $W_0^{1,2}(0, 1)$.

Throughout the paper, $C$ denotes a generic positive constant which may not be the same in each occurrence. More specific constants which depend on certain parameters $a,b,\dots$, are numbered as $C_{a,b,\dots}$. Unless otherwise stated, these constants may depend on the final time $T$, but remain independent of the discretization parameters $N$ and $\tau$, as well as the initial distribution.


\subsubsection{Well-posedness of SBE}

For the Dirichlet Laplacian $A$ and the Neumann Laplacian $A_1$, the operators $-A$ and $-A_1$ have eigenvalues \(\lambda_k=\pi^2k^2\) indexed by \(k\in\mathbb N_+\) and \(k\in\mathbb N\), respectively. The corresponding eigenfunctions are $e_k = \sqrt{2}\sin(k\pi\cdot)$ ($k \in \mathbb{N}_+$) for $-A$, and $h_0 = 1, h_k = \sqrt{2}\cos(k\pi\cdot)$ ($k \in \mathbb{N}_+$) for $-A_1$. The families \(\{e_k\}_{k\in\mathbb N_+}\) and \(\{h_k\}_{k\in\mathbb N}\) form complete orthonormal bases of $H$.

Define $\Lambda := I - A$ and $\Lambda_1 := I - A_1$. Both \(\Lambda\) and \(\Lambda_1\) are positive definite and self-adjoint operators with compact inverses. Their fractional powers can be defined for every $s \ge 0$ by
\begin{align}
  &\Lambda^{\frac{s}{2}} u= \sum_{k\in\N_+} (1+\lambda_k)^{\frac{s}{2}}\langle u,e_k\rangle_H e_k, \quad u\in \ID(\Lambda^{\frac{s}{2}}):=\{u\in H\colon \|u\|_{\HH^s}^2:=\sum_{k\in\N_+}(1+\lambda_k)^{s}\langle u,e_k\rangle_H^2<\infty\};\label{spectraldef:1}\\
   & \Lambda_1^{\frac{s}{2}} v= \sum_{k\in\N} (1+\lambda_k)^{\frac{s}{2}}\langle v,h_k\rangle_H h_k, \quad v\in \ID(\Lambda_1^{\frac{s}{2}}):=\{v\in H\colon \|v\|_{\HH_1^s}^2:=\sum_{k\in\N}(1+\lambda_k)^{s}\langle v,h_k\rangle_H^2<\infty\}.\label{spectraldef:2}
\end{align}
We set $\mathcal{H}^s:= \mathcal{D}(\Lambda^{\frac{s}{2}})$ and $\mathcal{H}_1^s := \mathcal{D}(\Lambda_1^{\frac{s}{2}})$
which form two families of separable Hilbert spaces. In particular, $\mathcal{H}_1^1 = W^{1,2}(0,1)$ and $\mathcal{H}_1^2 = \mathcal{D}(A_1) = \{v \in W^{2,2}(0,1)\colon \partial_x v(0) = \partial_x v(1) = 0\}$. 

Operators \(A\) and \(A_1\) generate the strongly continuous semigroups
$\{S(t):=e^{tA}\}_{t\ge0}$ and $\{S_1(t):=e^{tA_1}\}_{t\ge0}$ on \(H\), respectively. Their spectral representations are given by
\begin{align*}
	S(t)u=\sum_{k\in\N_+}e^{-t\lambda_k}\langle u, e_k\rangle_H e_k,\qquad S_1(t)u=\sum_{k\in\N}e^{-t\lambda_k}\langle u, h_k\rangle_H h_k,\qquad t\ge 0,~u\in H.
\end{align*}



Let $(\Omega,\mathcal F,\{\mathcal F_t\}_{t\in[0,T]},\mathbb P)$ be a complete filtered probability space satisfying the usual conditions. We write 
\begin{align}\label{KL-expansion}
W(t)=\sum_{k\in\N_+}\beta_k(t) e_k,\quad t\in[0,T]
\end{align}
formally for the cylindrical Wiener process on \(H\), where
\(\{\beta_k\}_{k\in\mathbb N_+}\) is a sequence of independent standard 
Wiener processes adapted to the filtration, with increments on $(s,t]$ independent of
$\mathcal F_s$ for $0\le s<t\le T$.  

For an \(U\)-valued random variable \(X\) and every \(p\in[1,\infty)\), define the  conditional norms
\[
\mathbb E_0[\|X\|_U^p]
:=\mathbb E[\|X\|_U^p|\mathcal F_0],
\qquad
\|X\|_{L_0^p(U)}
:=\bigl(\mathbb E_0[\|X\|_U^p]\bigr)^{1/p}.
\]
 For the SBE driven by the space--time white noise, extensive well-posedness results are available; see, e.g., \cite{daprato1994stochastic,gyongy1998existence}. The solution of \eqref{eq:sharp-sbe} is characterized in Lemma~\ref{lem:Burgers-rough-wellposed} and is used throughout the remainder of this paper.

\subsection{The Cole--Hopf transformation}\label{sec:2.2}

The classical Cole--Hopf substitution motivates the following
construction; see \cite{bertini1994stochastic,bertini1995stochastic}. Define
\begin{equation}\label{eq:sharp-canonical}
 U_0(x):=\exp\Bigl(\int_0^xu_0(y)\,dy
       -\min_{z\in[0,1]}\int_0^zu_0(y)\,dy\Bigr).
\end{equation}
Since $u_0\in H$ a.s., we have $U_0\in\HH_1^1$, $
\min_{x\in[0,1]}U_0(x)=1$, and $\partial_x\log U_0=u_0$ a.s. For \(\alpha\in[0,\frac14)\) and $k\in\N_+$, set
\[
b_{k,\alpha}(x):=\lambda_k^\alpha\int_0^x e_k(y)\,dy
=\sqrt{2}\pi^{2\alpha-1}k^{2\alpha-1}\bigl(1-\cos(k\pi x)\bigr),\quad x\in[0,1],
\]
and define the spatially integrated noise by $\widetilde W_\alpha(t,x)
:=
\sum_{k\in\mathbb N_+}
b_{k,\alpha}(x)\beta_k(t)$ for $t\in[0,T]$ and $x\in[0,1]$. It is a centered Gaussian field with covariance
\[
\E\Bigl[
\widetilde W_\alpha(t,x)
\widetilde W_\alpha(s,y)
\Bigr]
=
(t\wedge s)\sum_{k\in\mathbb N_+}
b_{k,\alpha}(x)b_{k,\alpha}(y).
\]
Set $\rho_\alpha(x):=
\frac12\sum_{k\in\mathbb N_+}b_{k,\alpha}(x)^2$ for $x\in[0,1]$. Then $d[\widetilde W_\alpha(\cdot,x)]_t
=
2\rho_\alpha(x)\,dt$. In particular, for the space-time white noise case, i.e., \(\alpha=0\), we have $\sum_{k\in\N_+}b_{k,0}(x)b_{k,0}(y)=x\wedge y$ and $\rho_0(x)=x/2$ for $x\in[0,1]$. 

A function $g$ satisfying $\partial_xg=u$ has the form $g(t,x)=\int_0^x u(t,y)\,dy+c(t).$ Formally, we can choose $c$ such that
\begin{equation}\label{eq:KPZ}
dg=\bigl(\partial_{xx}g+|\partial_xg|^2\bigr)\,dt
   +d\widetilde W_\alpha,
\qquad g(0,x):=\int_0^{x} u_0(y)\,dy
-\min_{z\in[0,1]}\int_0^z u_0(y)\,dy,
\end{equation}
with homogeneous Neumann boundary conditions. Thus, for $U=e^g$, \eqref{eq:KPZ} formally gives 
\[
dU
=U\,dg+\frac12U\,d[g]
=\left(\partial_{xx}U+\rho_\alpha U\right)dt
 +U\,d\widetilde W_\alpha.
\]
These calculations serve only as motivation. The rigorous
construction is given below.

\subsubsection{The linear heat equation with the multiplicative noise}

Set $U(t,x)=\exp(g(t,x))$ following the Cole--Hopf transformation in the preceding subsection. Formally, $U(t,x)$ satisfies the following stochastic heat equation with the linear multiplicative noise in the Stratonovich sense
\begin{equation}\label{eq:HC-transform}
\left\{
 \begin{aligned}
 &\partial_tU(t,x)=\partial_{xx} U(t,x)+U(t,x)\circ \partial_t\widetilde W_\alpha(t,x),\quad t\in(0,T], ~x\in(0,1),\\
 &U(0,\cdot)=U_0(\cdot):=\exp\Big(\int_0^{\cdot}u_0(y)\,dy-\min_{z\in[0,1]}\int_0^zu_0(y)\,dy\Big)=\exp(g(0,\cdot)).
\end{aligned}
\right.
\end{equation}
We interpret equation \eqref{eq:HC-transform} as the following It\^o-type stochastic evolution equation
\begin{equation}\label{eq:sharp-pam}
 dU=(A_1U+\rho_\alpha U)\,dt+G_\alpha(U)\,d\beta,
       \qquad U(0)=U_0,
\end{equation}
where $\rho_\alpha U$ represents the Stratonovich-to-It\^o correction. The stochastic integral
is defined as follows
\[
\beta(t)=\sum_{k\in\mathbb N_+}l_k\beta_k(t),\quad[G_\alpha(U)l_k](x) := U(x)b_{k,\alpha}(x), \quad U \in H,~k\in\N_+,
\]
where \(\{l_k\}_{k\in\mathbb N_+}\) denotes the canonical orthonormal basis of
\(\ell^2\).

The well-posedness of \eqref{eq:sharp-pam} and the analytical properties required for the inverse Cole--Hopf transformation and the subsequent error analysis are established in Section \ref{sec:linear}.


\subsection{Numerical method and main result}

 We now design a fully discrete approximation for the SBE \eqref{eq:sharp-sbe} by developing a robust scheme that leverages the Cole--Hopf transformation and an auxiliary truncation strategy. The construction consists of two steps: we first discretize the linear stochastic equation \eqref{eq:sharp-pam} for $U$, and then recover a numerical approximation of $u$ via a regularized inverse Cole--Hopf transformation.

We first introduce the spectral Galerkin method to obtain the spatial semi-discrete solution $\{U^N(t)\}_{t\in[0,T]}$ with discretization parameter $N$. For arbitrary $N\in \N_+$,  let $P_N$ and $\Pi_N$ be the orthogonal projection operators from $H$ onto $H_N:=\operatorname{span}\{e_1,\ldots,e_N\}$ and $H_{1,N}:=\operatorname{span}\{h_0,\ldots,h_N\}$, respectively, defined by
\[
P_Nu:=\sum_{j=1}^N\langle u,e_j\rangle_H e_j,\qquad \Pi_Nu:=\sum_{j=0}^N\langle u,h_j\rangle_H h_j,\qquad u\in H.
\]
For every $s\ge0$, $\Pi_N$ commutes with  \(\Lambda_1^{s/2}\) and is contractive on \(\mathcal H_1^s\). Moreover, for $0\le r\le s$, 
\begin{align}\label{ineq:projection}
 \|\Pi_Nu\|_{\HH_1^s}\le\|u\|_{\HH_1^s},\qquad
 \|(\mathrm I-\Pi_N)u\|_{\HH_1^r}
 \le (1+\lambda_{N+1})^{(r-s)/2}\|u\|_{\HH_1^s}.
\end{align}
 In particular, for all $\alpha_0\in[0,\infty)$, all  $u\in\HH_1^{\alpha_0}$, and all $N\in\N$, $\|(\mathrm I-\Pi_N)u\|_{H}\le (1+\lambda_{N+1})^{-\frac{\alpha_0}{2}}\|u\|_{\HH_1^{\alpha_0}}.$

The spectral Galerkin method applied to \eqref{eq:sharp-pam} is defined as the \(H_{1,N}\)-valued process \(U^N\) satisfying 
\begin{equation}\label{eq:Galerkin-pam}
\left\{
 \begin{aligned}
 &dU^N(t)=A_1U^N(t)\,dt+\Pi_N (\rho_\alpha U^N(t))\,dt+\Pi_N G_\alpha(U^N(t))\,d\beta(t),\quad t\in(0,T], \\
 &U^N(0)=\Pi_NU_0.
\end{aligned}
\right.
\end{equation}
For each $N$, the coefficients of this finite-dimensional stochastic differential equation are globally Lipschitz. Indeed, we have $\|\Pi_N(\rho_\alpha v)-\Pi_N(\rho_\alpha w)\|_H
\le
\|\rho_\alpha\|_\infty\|v-w\|_H$ and 
\[
 \|\Pi_NG_\alpha(v)-\Pi_NG_\alpha(w)\|_{\mathcal L_2(\ell^2;H)}^2=\sum_{k\in\N_+}\|\Pi_N((v-w)b_{k,\alpha})\|_H^2
 \le 2\int_0^1
\rho_\alpha(x)|v(x)-w(x)|^2\,dx\le
2\|\rho_\alpha\|_\infty
\|v-w\|_H^2.
\]
Using standard arguments (see, e.g., \cite{daprato2014stochastic}) gives a unique global \(H_{1,N}\)-valued mild solution
$\{U^N(t)\}_{t\in[0,T]}$, whose mild formulation is
\color{black}
\begin{align}\label{solution:HC-transform-SEE-spectralGalerkin}
U^N(t)\!=\!S_1(t)\Pi_NU_0+\int_0^t\!S_1(t-s)\Pi_N(\rho_\alpha U^N(s))\,ds+\sum_{k\in\N_+}\!\!\int_0^t\!S_1(t-s)\Pi_N(U^N(s)b_{k,\alpha})\,d\beta_k(s).
\end{align}

Fix $T>0$ and $N,M\ge 1$. In the temporal direction, an accelerated exponential Euler method (cf. \cite{Kruse2014Strong}) is further employed to obtain the fully discrete approximation $\{U_m^N\}_{m\in\{0,1,\ldots,M\}}$ with time step size $\tau:=T/M$ and $t_m:=m\tau$. We introduce the full discretization $\{U^N_m\}_{m\in\{0,1,\ldots,M\}}$ satisfying $U^N_0=\Pi_NU_0$ and 
\begin{align}\label{eq:sharp-update}
   U^N_{m}=S_1(\tau)U^N_{m-1}+ \int_{t_{m-1}}^{t_m}S_1(t_m-s)\Pi_N(\rho_\alpha U^N_{m-1})\,ds+\sum_{k\in\N_+}\!\!\int_{t_{m-1}}^{t_m}S_1(t_m-s)\Pi_N(U^N_{m-1}b_{k,\alpha})\,d\beta_k(s).
\end{align}

We next reconstruct a numerical approximation to the original SBE \eqref{eq:sharp-sbe} through a regularized inverse Cole--Hopf transformation. For constant $\eta>0$, define the mapping 
\begin{align*}
   \mathcal{R}_{N,\eta}(f):=P_N\left(\frac{\partial_x f}{f\vee\eta}\right),
   \qquad f\in \HH^1_1,
\end{align*}
where the maximum in the denominator is understood pointwise. For $\delta\in(0,1]$, define the mapping
\begin{equation}\label{eq:sharp-relative}
 \mathscr R_{N,\delta}(f)
  :=P_N\Big(\frac{\partial_x f}{f\vee(\delta\|f\|_{\HH_1^1})}\Big)\quad \mbox{for}~ f\ne0, \qquad \mbox{and} \qquad \mathscr R_{N,\delta}(0)=0.
\end{equation}
Equivalently,  
$\mathscr R_{N,\delta}(f)
=\mathcal R_{N,\delta}(f/\|f\|_{\HH_1^1})$ for $f\ne0$.
In particular, since 
\(\|\partial_xf\|_H\le\|f\|_{\mathcal H_1^1}\), we have
\begin{equation}\label{eq:sharp-cap}
 \mathscr R_{N,\delta}(cf)=\mathscr R_{N,\delta}(f)\quad\text{for constant}~ c>0,\qquad \mbox{and}
 \qquad \|\mathscr R_{N,\delta}(f)\|_H\le\delta^{-1}.
\end{equation}

 Fix $0<\varsigma<1/2-2\alpha$ and set $\delta:=1\wedge(N^{-\varsigma}+\sqrt\tau)$. The full discretization of the solution $\{u(t)\}_{t\in[0,T]}$ is then given by $u^{N,\delta}_0:=\mathcal{R}_{N,N^{-\varsigma}}(U^N_0)$ and
\begin{align}\label{eq:sharp-output}
u^{N,\delta}_m:=\mathscr{R}_{N,\delta}(U^N_m),\qquad m=1,2,\ldots,M.  
\end{align}

\subsubsection{Main result}

For a nonnegative random variable $X$ and for $p,r>0$, define
\begin{equation}\label{eq:sharp-budgets}
 \|X\|_{p,\infty}^p:=\sup_{R\ge1}R^p\mathbb P(X>R),\qquad \mathcal J_r:=\sup_{N\ge1}N^{2r}
    \E\big[\|u_0-\mathcal R_{N,N^{-r}}(U^N_0)\|_H^2\big].  
\end{equation}

Our main result gives the mean-square convergence rates of the fully discrete approximation for SBE~\eqref{eq:sharp-sbe}.

\begin{theo}\label{thm:sharp-main}
For arbitrary $T>0$ and $0\le \alpha<1/4$, let $\varsigma\in(0,1/2-2\alpha)$. Suppose that 
\[
\mathcal{J}_{\varsigma} + \big\|\|U_0\|_{\mathcal{H}_1^1}\big\|_{4,\infty}^4 <\infty. 
\]
 Then there exists a constant $C_{\alpha,\varsigma,T}>0$ such that, for all $N\in\N_+$ and all $\tau\in(0,1)$, the exact solution $u$ of \eqref{eq:sharp-sbe} with an $\mathcal{F}_0$-measurable initial value $u_0$ and the numerical approximation $\{u_m^{N,\delta}\}_{m\in\{0,1,\ldots,M\}}$ given by \eqref{eq:sharp-output} satisfy the error bound
\begin{equation*}
 \max_{0\le m\le M} \mathbb{E}\big[\|u(t_m) - u_m^{N,\delta}\|_H^2\big] 
 \le C_{\alpha,\varsigma,T}\big(1 + \mathcal{J}_\varsigma + \big\|\|U_0\|_{\mathcal{H}_1^1}\big\|_{4,\infty}^4\big)
\big( N^{-2\varsigma} + \tau(1 + 0\vee \ln(1/\tau))\big).
\end{equation*}
Here, $U_0$ is the Cole--Hopf transform of $u_0$ defined by~\eqref{eq:sharp-canonical}. 
\end{theo}

\begin{rem}\label{rem:sharp-Galerkin}
(1) For $0 < \alpha < 1/4$, the spatial convergence rate above is sharp within the spectral Galerkin approach. Indeed, the projection error of the stochastic convolution in the mild formulation of \eqref{eq:sharp-sbe} satisfies
\begin{equation*}
    \Big\|(I - P_N)\int_0^t S(t-s) (-A)^\alpha d W(s)\Big\|_{L^2(\Omega;H)} \asymp_{\alpha,t} N^{-1/2+2\alpha},\quad \forall t>0.
\end{equation*}
Therefore, the exponent $\varsigma<\frac12-2\alpha$ is the sharp spectral-Galerkin barrier. Moreover, $\alpha = 1/4$ represents a critical threshold: at $\alpha = 1/4$, the stochastic 
convolution in \eqref{eq:sharp-sbe} is not $L^2$-valued, while in the 
Cole--Hopf equation \eqref{eq:sharp-pam} the Stratonovich-to-It\^o correction 
$\rho_\alpha(x)\sim \sum_{k\in\mathbb{N}_+} k^{4\alpha-2}$ diverges logarithmically, 
both requiring renormalization beyond $\alpha = 1/4$. For the temporal discretization, the transformed solution has the temporal H\"older regularity of any order below \(1/2\) under the assumptions of the paper. Our accelerated scheme achieves a temporal convergence rate of order \(1/2\), up to a logarithmic correction, consistent with this regularity threshold. 

(2) Fix $0<\varsigma<1/2-2\alpha$, $\tau\in(0,1)$ and a initial law, there exists a constant $C_*\in(0,\infty)$ which may depend on $\alpha$, $\varsigma$, $T$, and the initial law but is independent of $N$ and $\tau$, such that the numerical scheme \eqref{eq:sharp-output} satisfies
\begin{equation*}
 \max_{0\le m\le M}\E\big[\|u(t_m)- u_m^{N,\delta}\|_H^2\big]
       \le C_*\left[N^{-2\varsigma}+\tau(1+0\vee \ln(1/\tau))\right].
\end{equation*}
Then, we have
\begin{equation*}
 \mathcal J_\varsigma\le C_*\qquad\mbox{and}\qquad\big\|\|u_0\|_{H}\big\|_{4,\infty}^4\le C_\varsigma(1+C_*).
\end{equation*}
In particular, $\mathcal{J}_\varsigma < \infty$ is strictly necessary for the error 
bound in Theorem \ref{thm:sharp-main}, and the weak-$L^4$ tail exponent $4$ is intrinsic for $L^2(\Omega;H)$ convergence. Moreover, since 
$\|u_0\|_H \le \|U_0\|_{\mathcal{H}_1^1} \le e^{2\|u_0\|_1}(1+\|u_0\|_H)$, the tail conditions on $\|U_0\|_{\mathcal{H}_1^1}$ and $\|u_0\|_H$ differ only through the exponential Cole--Hopf factor 
and coincide whenever $\int_0^\cdot u_0(y)\,dy$ is uniformly bounded.

(3) The logarithmic factor in Theorem \ref{thm:sharp-main} stems 
solely from the weak-$L^4$ tail truncation \eqref{ineq:tail1}. Whenever $\big\|\|U_0\|_{\mathcal{H}_1^1}\big\|_{4,\infty}^4 <\infty$ is strengthened to $\mathbb{E}[\|U_0\|_{\mathcal{H}_1^1}^4] < \infty$, the logarithmic factor can be removed.
\end{rem}

The proof of Remark~\ref{rem:sharp-Galerkin} (2) is deferred until after the proof of Theorem~\ref{thm:sharp-main}.

\section{Estimates for the transformed equation}\label{sec:linear}

This section establishes the analytical estimates for the transformed equation and its numerical error that will be used in the subsequent convergence analysis. In particular, we establish the moment bounds, strict positivity and inverse-moment estimates, together with the stability and regularity properties required for the inverse Cole--Hopf transformation.
\begin{lemma}\label{lem:sharp-existence}
Equation \eqref{eq:sharp-pam} has a unique mild solution $U$
with paths in $C([0,T];\HH_1^1)$. It is strictly positive on
$[0,T]\times[0,1]$ almost surely, and $u=\partial_x U/U$ is the unique $C([0,T];H)$ mild solution of \eqref{eq:sharp-sbe}. For every finite $p\ge2$ and $\gamma\in\{3/4,1\}$, the following estimates hold:
\begin{align}
 \sup_{0\le t\le T}\|U(t)\|_{L_0^p(\HH_1^\gamma)}
  +\sup_{0\le t\le T}\|U^N(t)\|_{L_0^p(\HH_1^\gamma)}
  +\max_{0\le m\le M}\|U_m^N\|_{L_0^p(\HH_1^\gamma)}&\le C_p\|U_0\|_{\HH_1^\gamma},
       \label{eq:sharp-moments}\\
 \Big\|\Big(\inf_{0\le t\le T,\,x\in[0,1]}U(t,x)\Big)^{-1}
          \Big\|_{L_0^p}&\le C_p,
       \label{eq:sharp-inverse-moments}\\
 \sup_{0\le t\le T}\|u(t)\|_{L_0^p(H)}&\le C_p\|U_0\|_{\HH_1^1}.
       \label{eq:sharp-u-moments}
\end{align}
\end{lemma}

\begin{proof}
Fix $r\in(0,1/2-2\alpha)$ and let
$\kappa:=1+r-\gamma\in(0,1)$.
We construct a solution of the mild equation
\begin{equation}\label{U mild}
U(t)=S_1(t) U_0+\int_0^t S_1(t-s) \rho_\alpha U(s) d s+\int_0^t S_1(t-s) G_\alpha(U(s)) d \beta(s),\qquad t\in[0,T].
\end{equation}
 Define
$F_{\gamma}(t):=\sup _{s\in[0,t]}\|U(s)\|_{L_0^p\left(\mathcal{H}_1^\gamma\right)}^2$ for $t\in[0,T]$. We deal with the convolutions in \eqref{U mild}.
Firstly, the conditional Burkholder--David--Gundy (BDG) inequality, Lemmas \ref{lem:analytic} and \ref{lem:sharp-multiplier} give
\begin{equation}\label{sto conv}
    \begin{aligned}
        & \left\|\int_0^t S_1(t-s) G_\alpha(U(s)) d \beta(s)\right\|_{L_0^p\left(\mathcal{H}_1^\gamma\right)}^2  \quad \leq C_p \int_0^t(t-s)^{\kappa-1}\left\|G_\alpha(U(s))\right\|_{L_0^p\left(\mathcal{L}_2\left(\ell^2 ; \mathcal{H}_1^r\right)\right)}^2 d s \\
& \quad \leq C_p \int_0^t(t-s)^{\kappa-1}\|U(s)\|_{L_0^p\left(\mathcal{H}_1^\gamma\right)}^2ds \quad \leq C_p\int_0^t(t-s)^{\kappa-1}F_{\gamma}(s)ds.
    \end{aligned}
\end{equation}
Second, for the deterministic convolution, we have
\begin{equation}\label{det conv}
    \begin{aligned}
& \left\|\int_0^t S_1(t-s) \rho_\alpha U(s) d s\right\|_{L_0^p\left(\mathcal{H}_1^\gamma\right)} \leq C \int_0^t(t-s)^{(\kappa-1)/2}\|U(s)\|_{L_0^p\left(\mathcal{H}_1^\gamma\right)}ds\leq C\int_0^t(t-s)^{(\kappa-1)/2} F_{\gamma}(s)^{1/2}ds.
    \end{aligned}
\end{equation}
Substituting \eqref{sto conv} and \eqref{det conv} into \eqref{U mild} and using $\left\|S_1(t) U_0\right\|_{\mathcal{H}_1^\gamma} \leq \|U_0\|_{\HH_1^\gamma},$ we obtain
\[
 F_{\gamma}(t)\le C\|U_0\|_{\HH_1^\gamma}^2+C\int_0^t(t-s)^{\kappa-1}F_{\gamma}(s)\,ds.
\]
Then it follows from Lemma \ref{lem:volterra} that
$$
\begin{aligned}
\sup _{0\leq s \leq t}\|U(s)\|_{L_0^p\left(\mathcal{H}_1^\gamma\right)}^2=F_{\gamma}(t) & \leq C_p \|U_0\|_{\HH_1^\gamma}^2+C_p \|U_0\|_{\HH_1^\gamma}^2 \int_0^t(t-s)^{\kappa-1} d s   \leq C_p \|U_0\|_{\HH_1^\gamma}^2 .
\end{aligned}
$$
Choose $t_0$ small enough such that $\frac{C_p}{\kappa}\left(t_0^\kappa+t_0^{\kappa+1}\right)<1.$ For two mild solutions $U$ and $V$
$$
\begin{aligned}
\sup _{t \leq t_0} \|U(t)-V(t)\|_{L_0^p\left(\mathcal{H}_1^\gamma\right)}^2 &\leq C_p\left[\int_0^{t_0}(t_0-s)^{\kappa-1} d s+t_0 \int_0^{t_0}(t_0-s)^{\kappa-1} d s\right] \sup _{s \leq t_0}\|U(s)-V(s)\|_{L_0^p\left(\mathcal{H}_1^\gamma\right)}^2 \\
& =\frac{C_p}{\kappa}\left(t_0^\kappa+t_0^{\kappa+1}\right) \sup _{s \leq t_0}\|U(s)-V(s)\|_{L_0^p\left(\mathcal{H}_1^\gamma\right)}^2 .
\end{aligned}
$$
 Iterating this estimate over consecutive small intervals gives the uniqueness of the mild solutions of \eqref{eq:sharp-pam}.

Similarly, the same convolution estimates give
\begin{equation}
    \begin{aligned}
        &\sup _{t \leq T}\left\|U^N(t)\right\|_{L_0^p\left(\mathcal{H}_1^\gamma\right)} \leq C_p \|U_0\|_{\HH_1^\gamma},\\
        &\left\|U_m^N\right\|_{L_0^p\left(\mathcal{H}_1^\gamma\right)}^2 \leq C_p \|U_0\|_{\HH_1^\gamma}^2+C_p \int_0^{t_m}\left(t_m-s\right)^{\kappa-1} \left\|U_{\lfloor s/\tau\rfloor}^N\right\|_{L_0^p\left(\mathcal{H}_1^\gamma\right)}^2ds,
    \end{aligned}
\end{equation}
where $\lfloor \cdot \rfloor$ denotes the floor function: if $m\ge 0$ and $t\in[t_m,t_{m+1})$, one has $\lfloor t/\tau \rfloor=m$. Thus \eqref{eq:sharp-moments} is proved by setting $E_m=\left\|U_m^N\right\|_{L_0^p\left(\mathcal{H}_1^\gamma\right)}^2$ in Lemma \ref{lem:volterra}.

To show $U(\cdot)\in C([0,T],\mathcal{H}_1^1)$, let $0<\eta<\frac{r}{2},~ p>\frac{1}{\eta},$ and denote the stochastic convolution in \eqref{U mild} by
$$
Y_\eta(t):=\int_0^t(t-s)^{-\eta} S_1(t-s) G_\alpha(U(s)) d \beta(s),
$$
where using \eqref{eq:sharp-moments} and Lemma \ref{lem:sharp-multiplier}, we obtain
\begin{equation}
    \begin{aligned}
        \sup _{s \leq T}\|G_\alpha(U(s))\|_{L_0^p\left(\mathcal{L}_2\left(\ell^2 ; \mathcal{H}_1^r\right)\right)} \leq C  \sup _{s \leq T}\|U(s)\|_{L_0^p\left(\mathcal{H}_1^1\right)} \leq C_{p}\left\|U_0\right\|_{\mathcal{H}_1^1}.
    \end{aligned}
\end{equation}
By Lemma \ref{lem:analytic}
$$
\begin{aligned}
(t-s)^{-2 \eta}\left\|S_1(t-s)\right\|_{\mathcal{L}\left(\mathcal{H}_1^r ; \mathcal{H}_1^1\right)}^2 & \leq C(t-s)^{-2 \eta-(1-r)}  =C(t-s)^{r-1-2 \eta} .
\end{aligned}
$$
The conditional BDG inequality gives
$$
\begin{aligned}
\left\|Y_\eta(t)\right\|_{L_0^p\left(\mathcal{H}_1^1\right)}^2 & \leq C_p\|U_0\|_{\mathcal{H}_1^1}^2  \int_0^t(t-s)^{r-1-2 \eta} d s  =\frac{C_p }{r-2 \eta} t^{r-2 \eta}\|U_0\|_{\mathcal{H}_1^1}^2 .
\end{aligned}
$$
The factorization method
\cite{daprato1988regularity}  yields
\[
 \int_0^tS_1(t-s)G_\alpha(U(s))\,d\beta(s)
  =\frac{\sin(\pi\eta)}{\pi}
       \int_0^t(t-s)^{\eta-1}S_1(t-s)Y_\eta(s)\,ds.
\]
Notice that $p>\frac{1}{\eta}  \Longleftrightarrow (\eta-1) \frac{p}{p-1}>-1$, Hölder's inequality gives
$$
\begin{aligned}
& \left\|\int_0^{\cdot}(\cdot-s)^{\eta-1} S_1(\cdot-s) Y_\eta(s) d s\right\|_{C\left([0, T] ; \mathcal{H}_1^1\right)}  \quad \leq C_{\eta, p, T}\left(\int_0^T\left\|Y_\eta(s)\right\|_{\mathcal{H}_1^1}^p d s\right)^{1 / p}.
\end{aligned}
$$
Taking conditional moments yields
$$
\left\|\int_0^\cdot S_1(\cdot-s) G_\alpha(U(s))d \beta(s)\right\|_{L_0^p\left(C\left([0, T] ; \mathcal{H}_1^1\right)\right)} \leq C_p\Vert{}U_0\Vert{}_{\mathcal{H}_1^1}.
$$
For the deterministic convolution in \eqref{U mild}, Lemma \ref{lem:analytic} gives $U \in C\left([0, T] ; \mathcal{H}_1^1\right)$ by the following
$$
\begin{aligned}
& \Big\|\int_0^t S_1(t-s) \rho_\alpha U(s) d s\Big\|_{\mathcal{H}_1^1}  \quad \leq C \int_0^t(t-s)^{(r-1) / 2}\left\|\rho_\alpha U(s)\right\|_{\mathcal{H}_1^r} ds \\
& \quad \leq C \sup _{s\in[0,T]}\|U(s)\|_{\mathcal{H}_1^r} \int_0^t(t-s)^{(r-1) / 2} d s  \quad=\frac{2 C}{r+1} t^{(r+1) / 2} \sup _{s\in[0,T]}\|U(s)\|_{\mathcal{H}_1^r} .
\end{aligned}
$$
For the strict positiveness of $U$, Lemma~\ref{lem:FeynmanKac} gives
\[
 \hat{U}^K(t,x)=\E_R\Bigl[U_0(R^x(t))
     \exp\Bigl\{\sum_{k\le K}\int_0^t
                    b_{k,\alpha}(R^x(t-s))\,d\beta_k(s)\Bigr\}\Bigr],
\]
where $R^x$ is the reflected Brownian motion on $[0,1]$ with generator $\partial_{xx}$,
independent of the noise. Define
$$
Z_K(t):=\sum_{k \leq K} \int_0^t S_1(t-s) b_{k,\alpha} d \beta_k(s).
$$
Using $U_0 \geq 1$ and Jensen's inequality gives
$$
\begin{aligned}
&\hat{U}^K(t, x)  \geq \exp \left\{\mathbb{E}_R\left[\sum_{k \leq K} \int_0^t b_{k,\alpha}\left(R^x(t-s)\right) d \beta_k(s)\right]\right\} \\
& =\exp \left\{\sum_{k \leq K} \int_0^t \mathbb{E}_R\left[b_{k,\alpha}\left(R^x(t-s)\right)\right] d \beta_k(s)\right\}  =e^{Z_K} .
\end{aligned}
$$
Similar factorization calculation yields
$Z_K \longrightarrow Z:=\sum_{k \geq 1} \int_0^\cdot S_1(\cdot-s) b_{k,\alpha} d \beta_k(s) \text { in } L_0^p\left(C\left([0, T] ; \mathcal{H}_1^1\right)\right)$ and similarly $\hat{U}^K \longrightarrow U$ in $L_0^p(C([0,T];\mathcal{H}_1^1))$. Sobolev embedding gives $\|f\|_{\infty} \leq C\|f\|_{\mathcal{H}_1^1}$. Therefore a subsequence satisfies
$$
\hat{U}^{K_j} \longrightarrow U, \quad Z_{K_j} \longrightarrow Z \quad \text { in } C([0, T] \times[0,1]),\qquad a.s.
$$
Consequently, $U(t, x) \geq e^{Z(t, x)}>0  \text { for every }(t, x) \in[0, T] \times[0,1] .$ Hence
$$
\left(\inf _{(t,x) \in[0, T] \times[0,1]} U(t,x)\right)^{-1} \leq e^{-\inf Z} \leq e^{\|Z\|_{C(0, T] \times(0,1)}} .
$$
Since $Z$ is Gaussian and independent of $\mathcal{F}_0$, Fernique's theorem gives $a>0$ such that $
\mathbb{E}\big[\exp \left(a\|Z\|_{\infty}^2\right)\big]<\infty.$
Thus \eqref{eq:sharp-inverse-moments} is proved by
$$
\begin{aligned}
\mathbb{E}_0\Big[\Big(\inf_{(t, x) \in[0, T] \times[0,1]} U(t, x)\Big)^{-p}\Big]\le \mathbb{E} [e^{p\|Z\|_{\infty}}] 
 \leq e^{p^2 /(4 a)} \mathbb{E}[ e^{a\|Z\|_{\infty}^2}] <\infty .
\end{aligned}
$$
Notice that
$$
\begin{aligned}
\|u(t)\|_H & =\left\|\frac{\partial_xU(t)}{U(t)}\right\|_H  \leq\left\|\partial_xU(t)\right\|_H\left(\inf _{x\in[0,1]} U(t, x)\right)^{-1}  \leq\|U(t)\|_{\mathcal{H}_1^1}\left(\inf _{(s,x) \in[0, T] \times[0,1]} U(s, x)\right)^{-1} .
\end{aligned}
$$
Using the H\"older inequality, \eqref{eq:sharp-moments} and \eqref{eq:sharp-inverse-moments} gives \eqref{eq:sharp-u-moments} by the following
$$
\begin{aligned}
\|u(t)\|_{L_0^p(H)} & \leq\|U(t)\|_{L_0^{2 p}\left(\mathcal{H}_1^1\right)}\left\|\left(\inf _{(t,x) \in[0, T] \times[0,1]} U(t, x)\right)^{-1}\right\|_{L_0^{2 p}} \leq C_p \|U_0\|_{\HH_1^1} .
\end{aligned}
$$
Finally, we justify $u=\partial_x U/U$. Since   $\rho_\alpha=\frac12\sum_kb_{k,\alpha}^2$ is bounded, for
$\varphi\in\HH_1^1$, Lemma~\ref{lem:finite-noise-log-chain} gives
\begin{equation}\label{eq:sharp-log-weak}
 \begin{split}
 d\langle\log U,\varphi\rangle_H
 &=\left[-\langle\partial_x\log U,\partial_x\varphi\rangle_H
      +\int_0^1|\partial_x\log U|^2\varphi\,dx\right]dt
      +\sum_k\langle b_{k,\alpha},\varphi\rangle_H\,d\beta_k.
 \end{split}
\end{equation}
Testing \eqref{eq:sharp-log-weak} against $h_k$ and using
$\langle b_{j,\alpha}, h_k \rangle_H = -\delta_{jk}\lambda_k^\alpha / (k\pi)$
give the mild equation
\begin{equation}\label{eq:sharp-mild}
 u(t)=S(t)u_0+\int_0^tS(t-s)\partial_x(u(s)^2)\,ds
       +\sum_{k\ge1}\lambda_k^\alpha
                   \int_0^te^{-\lambda_k(t-s)}e_k\,d\beta_k(s).
\end{equation}
Using Lemma \ref{lem:Burgers-rough-wellposed} completes the proof.
\end{proof}

\begin{lemma}\label{lem:sharp-linear-error}
Let $\gamma\in\{3/4,1\}$ and $0<a<1/2-2\alpha$. For every finite $p\ge2$ and $t,t_m\le1\wedge T$
\begin{align}
 \|U(t)-S_1(t)U_0\|_{L_0^p(\HH_1^\gamma)}
       &\le C_p\|U_0\|_{\HH_1^\gamma} t^{a/2},
       \label{eq:sharp-start-exact}\\
 \|U_m^N-S_1(t_m)\Pi_NU_0\|_{L_0^p(\HH_1^\gamma)}
       &\le C_p\|U_0\|_{\HH_1^\gamma} t_m^{a/2},
       \label{eq:sharp-start-grid}\\
       \|U(t_m)-U_m^N\|_{L_0^p(\HH_1^\gamma)}
       &\le C_p\|U_0\|_{\HH_1^\gamma}
           (e^{-\lambda_{N+1}t_m}+N^{-a}+\sqrt\tau).
       \label{eq:sharp-linear-error}
\end{align}
\end{lemma}
\begin{proof}
Choose $a<r<1/2-2\alpha$ and let $\kappa:=1+r-\gamma\in(0,1)$. The mild equation to \eqref{eq:sharp-pam} is given by
$$
U(t)-S_1(t) U_0=\int_0^t S_1(t-s) \rho_\alpha U(s) d s+\int_0^t S_1(t-s) G_\alpha(U(s)) d \beta(s).
$$
Using Lemma \ref{lem:analytic} and  Lemma \ref{lem:sharp-existence} yields
$$
\begin{aligned}
&\Big\|\int_0^t S_1(t-s) \rho_\alpha U(s) d s\Big\|_{L_0^p\left(\mathcal{H}_1^\gamma\right)}  \leq C \int_0^t(t-s)^{(\kappa-1) / 2}\|U(s)\|_{L_0^p\left(\mathcal{H}_1^r\right)} d s \\
& \leq C_p \|U_0\|_{\HH_1^\gamma} \int_0^t(t-s)^{(\kappa-1) / 2} d s 
=\frac{ C_p}{\kappa+1} \|U_0\|_{\HH_1^\gamma} t^{(\kappa+1) / 2} .
\end{aligned}
$$
The conditional BDG inequality gives
$$
\begin{aligned}
&\Big\|\int_0^t S_1(t-s) G_\alpha(U(s)) d \beta(s)\Big\|_{L_0^p\left(\mathcal{H}_1^\gamma\right)}^2  \leq C_p \int_0^t(t-s)^{\kappa-1}\left\|G_\alpha(U(s))\right\|_{L_0^p\left(\mathcal{L}_2\left(\ell^2 ; \mathcal{H}_1^r\right)\right)}^2 d s \\
& \leq C_p \|U_0\|_{\HH_1^\gamma}^2 \int_0^t(t-s)^{\kappa-1} d s  =\frac{C_p}{\kappa} \|U_0\|_{\HH_1^\gamma}^2 t^\kappa .
\end{aligned}
$$
Consequently, for $\kappa>a$ and $0<t \leq 1 \wedge T$, it holds that
$$
\begin{aligned}
\left\|U(t)-S_1(t) U_0\right\|_{L_0^p\left(\mathcal{H}_1^\gamma\right)} & \leq C_p \|U_0\|_{\HH_1^\gamma}\left(t^{(\kappa+1) / 2}+t^{\kappa / 2}\right)  \leq C_p \|U_0\|_{\HH_1^\gamma} t^{\kappa / 2} \leq C_p\|U_0\|_{\HH_1^\gamma} t^{a / 2}. 
\end{aligned}
$$

For the discrete solution $U_m^N$,  we have
$$
\begin{aligned}
U_m^N-S_1\left(t_m\right) U^N_0= & \int_0^{t_m} S_1\left(t_m-s\right) \Pi_N \rho_\alpha U_{\lfloor s/\tau\rfloor}^N d s  +\int_0^{t_m} S_1\left(t_m-s\right) \Pi_N G_\alpha\left(U_{\lfloor s/\tau\rfloor}^N\right) d \beta(s) .
\end{aligned}
$$
Using $\|U_{\lfloor s/\tau\rfloor}^N\|_{L_0^p\left(\mathcal{H}_1^\gamma\right)} \leq C_p \|U_0\|_{\HH_1^\gamma}$ and repeating the calculations above yield
$$
\begin{aligned}
\left\|U_m^N-S_1\left(t_m\right) U^N_0\right\|_{L_0^p\left(\mathcal{H}_1^\gamma\right)} & \leq C_p \|U_0\|_{\HH_1^\gamma}\left(t_m^{(\kappa+1) / 2}+t_m^{\kappa / 2}\right)  \leq C_p \|U_0\|_{\HH_1^\gamma} t_m^{a / 2}, \quad t_m \leq 1 \wedge T,
\end{aligned}
$$
which proves \eqref{eq:sharp-start-grid}.

For \eqref{eq:sharp-linear-error}, we decompose $U(t_m)-U_m^N=U(t_m)-U^N(t_m)+U^N(t_m)-U_m^N,$ which represent the spatial error and temporal error correspondingly. For the spatial error, mild formulations give
\begin{equation}\label{UUN}
\begin{aligned}
    U(t)-U^N(t)&=S_1(t)\left(\mathrm{I}-\Pi_N\right) U_0+\int_0^t S_1(t-s) \Pi_N \rho_\alpha (U-U^N) \,d s\\
    &\quad +\int_0^t S_1(t-s)\left(\mathrm{I}-\Pi_N\right) \rho_\alpha U \,d s+\int_0^t S_1(t-s) \Pi_N G_\alpha(U-U^N) \,d \beta(s)\\
    &\quad +\int_0^t S_1(t-s)\left(\mathrm{I}-\Pi_N\right) G_\alpha(U) \,d \beta(s)\\
    &=:J_0+J_1+J_2+J_3+J_4.
    \end{aligned}
\end{equation}
For the term $J_0$, we derive
\begin{equation}\label{J0}
\left\|{J}_0(t)\right\|_{\mathcal{H}_1^\gamma}^2  =\sum_{k>N}\left(1\!+\!\lambda_k\right)^\gamma e^{-2 \lambda_k t}\left|\left\langle U_0, h_k\right\rangle_H\right|^2  
\leq e^{-2 \lambda_{N+1} t} \|U_0\|_{\HH_1^\gamma}^2 .
\end{equation}
For the term $J_1$,  Lemmas \ref{lem:sharp-multiplier} and  \ref{lem:analytic} yield
\begin{equation}\label{J1}
    \begin{aligned}
\left\|{J}_1(t)\right\|^2_{L_0^p\left(\mathcal{H}_1^\gamma\right)}=& \Big\| \int_0^t S_1(t-s) \Pi_N \rho_\alpha (U-U^N) d s \Big\|^2_{L_0^p\left(\mathcal{H}_1^\gamma\right)} \\
 \leq& C \Big(\int_0^t(t-s)^{(\kappa-1) / 2}\left\|U-U^N\right\|_{L_0^p\left(\mathcal{H}_1^\gamma\right)} d s\Big)^2 \leq C_T \int_0^t(t-s)^{\kappa-1}\left\|U-U^N\right\|^2_{L_0^p\left(\mathcal{H}_1^\gamma\right)} d s 
\end{aligned}
\end{equation}
Define $\theta:=\gamma+a-r<1$. For $k>N$, it holds that $\left(1+\lambda_k\right)^{\gamma-r} e^{-2 \lambda_k t} \leq C \lambda_k^{-a} t^{-\theta} \leq C N^{-2 a} t^{-\theta} .$ Therefore
\begin{equation}\label{theta ine}
\left\|\left(\mathrm{I}-\Pi_N\right) S_1(t)\right\|_{\mathcal{L}\left(\mathcal{H}_1^r ; \mathcal{H}_1^\gamma\right)} \leq C N^{-a} t^{-\theta / 2}.
\end{equation}
For the term $J_2$, then applying Lemma \ref{lem:sharp-multiplier}, \eqref{eq:sharp-moments} and \eqref{theta ine} gives
\begin{equation}\label{J2}
 \left\|{J}_2(t)\right\|_{L_0^p\left(\mathcal{H}_1^\gamma\right)}  \leq C N^{-a} \int_0^t(t-s)^{-\theta / 2}\left\|\rho_\alpha U\right\|_{L_0^p\left(\mathcal{H}_1^r\right)} d s  \leq C_p N^{-a} \|U_0\|_{\HH_1^\gamma} \int_0^t(t-s)^{-\theta / 2} d s .
\end{equation}
For the term $J_3$, by the conditional BDG inequality and Lemma \ref{lem:sharp-multiplier}, we obtain
\begin{equation}\label{J3}
\begin{aligned}
& \left\|{J}_3(t)\right\|_{L_0^p\left(\mathcal{H}_1^\gamma\right)}^2  \leq C_p\!\! \int_0^t\!(t\!-\!s)^{\kappa\!-\!1}\left\|\Pi_N G_\alpha\left(U\!-\!U^N\right)\right\|_{L_0^p\left(\mathcal{L}_2\left(\ell^2 ; \mathcal{H}_1^r\right)\right)}^2 d s \leq C_p\!\! \int_0^t\!(t\!-\!s)^{\kappa\!-\!1}\left\|U\!-\!U^N\right\|_{L_0^p\left(\mathcal{H}_1^\gamma\right)}^2 d s.
\end{aligned}
\end{equation}
Similarly,
\begin{equation}\label{J4}
\begin{aligned}
& \left\|{J}_4(t)\right\|_{L_0^p\left(\mathcal{H}_1^\gamma\right)}^2  \leq C_p N^{-2 a} \int_0^t(t-s)^{-\theta}\left\|G_\alpha(U)\right\|_{L_0^p\left(\mathcal{L}_2\left(\ell^2 ; \mathcal{H}_1^r\right)\right)}^2 d s  \leq C_p N^{-2 a} \|U_0\|_{\HH_1^\gamma}^2 \int_0^t(t-s)^{-\theta} d s  .
\end{aligned}
\end{equation}
Substituting \eqref{J0}-\eqref{J4} into \eqref{UUN} yields
$$
\left\|U(t)-U^N(t)\right\|_{L_0^p\left(\mathcal{H}_1^\gamma\right)}^2 \leq C_p \|U_0\|_{\HH_1^\gamma}^2\left(e^{-2 \lambda_{N+1} t}+N^{-2 a}\right)+C_p \int_0^t(t-s)^{\kappa-1} \left\|U(s)-U^N(s)\right\|_{L_0^p\left(\mathcal{H}_1^\gamma\right)}^2 d s.
$$
Lemma \ref{lem:volterra} gives
\begin{equation}\label{spa err}
\begin{aligned}
\left\|U(t)-U^N(t)\right\|_{L_0^p\left(\mathcal{H}_1^\gamma\right)}^2 \leq& C_p \|U_0\|_{\HH_1^\gamma}^2  {\Big[e^{-2 \lambda_{N+1} t}+N^{-2 a}} +\int_0^t(t-s)^{\kappa-1}\left(e^{-2 \lambda_{N+1} s}+N^{-2 a}\right) d s\Big]\\
\leq& C_p \|U_0\|_{\HH_1^\gamma}^2\left(e^{-2 \lambda_{N+1} t}+N^{-2 a}\right).
\end{aligned}
\end{equation}
Denote $\ell(s):=\tau\lfloor \frac{s}{\tau}\rfloor$ for $s\in[0,T)$, so that \(\ell(s)=t_j\) for \(s\in[t_j,t_{j+1})\). For the temporal error $U^N(t_m)-U^N_m$, we have
\begin{equation}\label{4K}
\begin{aligned}
&U^N(t_m)-U^N_m=  \int_0^{t_m} S_1\left(t_m-s\right) \Pi_N \rho_\alpha\left(U^N(s)-U_{\lfloor s/\tau\rfloor}^N\right) d s \\
&\qquad +\int_0^{t_m} S_1\left(t_m-s\right) \Pi_N G_\alpha\left(U^N(s)-U_{\lfloor s/\tau\rfloor}^N\right) d \beta(s)\\
&=\int_0^{t_m} S_1\left(t_m-s\right) \Pi_N \rho_\alpha \left(U^N(s)-U^N(\ell(s))\right) d s +\int_0^{t_m} S_1\left(t_m-s\right) \Pi_N \rho_\alpha \left( U^N(\ell(s))-U_{\lfloor s/\tau\rfloor}^N \right) d s \\
&\quad +\int_0^{t_m} S_1\left(t_m-s\right) \Pi_N G_\alpha\left(U^N(s)-U^N(\ell(s))\right)  d \beta(s)+\int_0^{t_m} S_1\left(t_m-s\right) \Pi_N G_\alpha\left( U^N(\ell(s))-U_{\lfloor s/\tau\rfloor}^N \right) d \beta(s)\\
&=:K_1+K_2+K_3+K_4.
\end{aligned}
\end{equation}
We first claim that
\begin{equation}\label{neq1}
\int_0^{t_m}\left(t_m-s\right)^{\kappa-1} \|U^N(s)-U^N(\ell(s)) \|^2_{L_0^p\left(\mathcal{H}_1^r\right)} d s \leq C_{p, T}\|U_0\|_{\HH_1^\gamma}^2 \tau .
\end{equation}
Indeed,
\[
U^N(s)-U^N(\ell(s))
=U^N(s)-S_1(s)U^N_0-U^N(\ell(s))+S_1(\ell(s))U^N_0
+\bigl(S_1(s)-S_1(\ell(s))\bigr)U^N_0 .
\]
Splitting the deterministic and stochastic convolutions at $\ell(s)$ gives
\[
\left\|U^N(s)-S_1(s)U^N_0-U^N(\ell(s))+S_1(\ell(s))U^N_0\right\|_{L_0^p(\HH_1^r)}^2
\le C_{p,T}\|U_0\|_{\HH_1^\gamma}^2(s-\ell(s)).
\]
Moreover, Lemma~\ref{lem:analytic} and
$r+1-\gamma=\kappa$ give, for $s\in[t_j,t_{j+1})$,
\[
\left\|\bigl(S_1(s)-S_1(\ell(s))\bigr)U^N_0\right\|_{\HH_1^r}^2
\le C\|U_0\|_{\HH_1^\gamma}^2
\begin{cases}
s^{1-\kappa},&j=0,\\
(s-t_j)t_j^{-\kappa},&j\ge1 .
\end{cases}
\]
Consequently
\[
\begin{aligned}
&\int_0^{t_m}(t_m-s)^{\kappa-1}
 \left\|U^N(s)-U^N(\ell(s))\right\|_{L_0^p(\HH_1^r)}^2ds\\
&\le C_{p,T}\|U_0\|_{\HH_1^\gamma}^2\Bigg[
 \int_0^{t_m}(t_m-s)^{\kappa-1}(s-\ell(s))\,ds
 +\int_0^\tau(t_m-s)^{\kappa-1}s^{1-\kappa}\,ds+\tau\sum_{j=1}^{m-1}\int_{t_j}^{t_{j+1}}
 (t_m-s)^{\kappa-1}t_j^{-\kappa}\,ds\Bigg].
\end{aligned}
\]
The three integrals are estimated as follows. First,
\[
\int_0^{t_m}(t_m-s)^{\kappa-1}(s-\ell(s))\,ds
\le \tau\int_0^{t_m}(t_m-s)^{\kappa-1}\,ds
=\frac{\tau t_m^\kappa}{\kappa}
\le C_{\kappa,T}\tau .
\]
Second,
\[
\int_0^\tau(t_m-s)^{\kappa-1}s^{1-\kappa}\,ds
\le
\begin{cases}
\displaystyle
\tau\int_0^1(1-y)^{\kappa-1}y^{1-\kappa}\,dy,&m=1,\\[6pt]
\displaystyle
\tau^{\kappa-1}\int_0^\tau s^{1-\kappa}\,ds
=\frac{\tau}{2-\kappa},&m\ge2,
\end{cases}
\le C_\kappa\tau .
\]
Finally, $s-t_j\le\tau$ and $t_j\ge s/2$ for $j\ge1$ imply
\[
\begin{aligned}
\tau\sum_{j=1}^{m-1}\int_{t_j}^{t_{j+1}}
 (t_m-s)^{\kappa-1}t_j^{-\kappa}\,ds
&\le C_\kappa\tau\int_0^{t_m}
 (t_m-s)^{\kappa-1}s^{-\kappa}\,ds=C_\kappa\tau\int_0^1(1-y)^{\kappa-1}y^{-\kappa}\,dy
\le C_\kappa\tau .
\end{aligned}
\]
Combining these three bounds proves \eqref{neq1}. For the term $K_1$, using Lemmas \ref{lem:sharp-multiplier} and \ref{lem:analytic}, and \eqref{neq1}, we have
\begin{equation}\label{K1}
\left\|{K}_{1}\right\|_{L_0^p\left(\mathcal{H}_1^\gamma\right)} \leq C \int_0^{t_m}\left(t_m-s\right)^{(\kappa-1) / 2}\left\|U^N(s)-U^N(\ell(s))\right\|_{L_0^p\left(\mathcal{H}_1^r\right)} d s \leq C_{p, T}\|U_0\|_{\HH_1^\gamma} \tau^{\frac12}.
\end{equation}
For the term $K_2$, direct computation gives
\begin{equation}\label{K2}
 \left\|{K}_{ 2}\right\|_{L_0^p\left(\mathcal{H}_1^\gamma\right)}  \leq C \int_0^{t_m}\left(t_m-s\right)^{(\kappa-1) / 2}\left\| U^N(\ell(s))-U_{ \lfloor s/\tau\rfloor}^N \right\|_{L_0^p\left(\mathcal{H}_1^\gamma\right)} d s.
\end{equation}
For the term $K_3$, the conditional BDG inequality and \eqref{neq1} give
\begin{equation}\label{K3}
\begin{aligned}
 \left\|{K}_{ 3}\right\|_{L_0^p\left(\mathcal{H}_1^\gamma\right)}^2  \leq& C_p \int_0^{t_m}\left(t_m-s\right)^{\kappa-1}\left\|\Pi_N G_\alpha\left(U^N(s)-U^N(\ell(s))\right)\right\|_{L_0^p\left(\mathcal{L}_2\left(\ell^2 ; \mathcal{H}_1^r\right)\right)}^2 d s \\
 \leq& C_p \int_0^{t_m}\left(t_m-s\right)^{\kappa-1}\left\|U^N(s)-U^N(\ell(s))\right\|_{L_0^p\left(\mathcal{H}_1^r\right)}^2 d s  \leq C_{p, T} \|U_0\|_{\HH_1^\gamma}^2 \tau.
\end{aligned}
\end{equation}
For the term $K_4$, it holds that
\begin{equation}\label{K4}
\begin{aligned}
\left\|{K}_{ 4}\right\|_{L_0^p\left(\mathcal{H}_1^\gamma\right)}^2 
& \leq C_p \int_0^{t_m}\left(t_m-s\right)^{\kappa-1}\left\|\Pi_N G_\alpha\left(U^N(\ell(s))-U_{\lfloor s/\tau\rfloor}^N \right)\right\|_{L_0^p\left(\mathcal{L}_2\left(\ell^2 ; \mathcal{H}_1^r\right)\right)}^2 d s \\
& \leq C_p \int_0^{t_m}\left(t_m-s\right)^{\kappa-1}\left\|U^N(\ell(s))-U_{\lfloor s/\tau\rfloor}^N \right\|_{L_0^p\left(\mathcal{H}_1^\gamma\right)}^2 d s.
\end{aligned}
\end{equation}
Substituting \eqref{K1}--\eqref{K4} into \eqref{4K}, we obtain
\begin{equation*}
    \|U^N(t_m)-U^N_m\|_{L_0^p\left(\mathcal{H}_1^\gamma\right)}^2 \leq C_{p, T} \|U_0\|_{\HH_1^\gamma}^2 \tau+C_{p, T} \int_0^{t_m}\left(t_m-s\right)^{\kappa-1} \left\|U^N(\ell(s))-U_{\lfloor s/ \tau\rfloor}^N \right\|_{L_0^p\left(\mathcal{H}_1^\gamma\right)}^2ds.
\end{equation*}
Using Lemma \ref{lem:volterra} into the above inequality with $K=C_{p, T} \|U_0\|_{\HH_1^\gamma}^2 \tau$ yields
\begin{equation}\label{tem err}
\max _{0 \leq m \leq M}\left\|U^N\left(t_m\right)-U_m^N\right\|_{L_0^p\left(\mathcal{H}_1^\gamma\right)} \leq C_{p, T} \|U_0\|_{\HH_1^\gamma} \sqrt{\tau} .
\end{equation}
Combining \eqref{spa err} and \eqref{tem err} proves \eqref{eq:sharp-linear-error}.
\end{proof}

\section{Convergence analysis}\label{sec:reconstruction}

We now prove the main convergence result given in Theorem \ref{thm:sharp-main}. The proof combines the error estimates for the transformed equation obtained in Section \ref{sec:linear} and a stability estimate for the regularized inverse Cole--Hopf transformation. A localization technique, together with the weak-tail control encoded by \(\big\|\|U_0\|_{\mathcal{H}_1^1}\big\|_{4,\infty}\), is then used to handle the random initial datum. 


\subsection{Auxiliary estimates}

The main difficulty in the convergence analysis is to transfer the approximation error of the transformed variables to that of the original Burgers solution. The inverse Cole--Hopf transformation $f\mapsto \partial_x f/f$ is nonlinear and, more importantly, its stability depends on the lower bound of the state variable \(f\). Hence, a small error in the transformed variable does not directly imply a small error after inversion when \(f\) approaches zero. The following lemma quantifies this error propagation by combining the approximation error with suitable inverse-moment control.

Let $\mathcal{V}\in\IF_0$ and let $f_1,f_2$ be $\HH_1^1$-valued random variables.  Assume that, a.s., on $\mathcal{V}$, the continuous modification of $f_2$ satisfies $f_2(x)>0$ for every $x\in[0,1]$. For a Banach space \(\mathbb{B}\) and a \(\mathbb{B}\)-valued random variable \(X\), define
$$
 \|X\|_{\mathbb{B},p,\mathcal{V}}:=(\E_0[\|X\|_{\mathbb{B}}^p\mathbf1_{\mathcal{V}}])^{1/p},\qquad p\ge 2,
$$
where $\|\cdot\|_{\mathbb{B}}$ denotes the norm on \(\mathbb{B}\).

\begin{lemma}\label{lem:error-f1f2}
Set $\Psi(z):=z^2(1+0\vee \log(1/z))$ for $z>0$ and $\Psi(0)=0$. Suppose that, for every $p\ge2$,
\begin{align}\label{conditions-f1f2}
 \|f_1\|_{\HH_1^1,p,\mathcal{V}}+\Big\|\Big\|\frac{\partial_x f_2}{f_2}\Big\|_H\Big\|_{\R,p,\mathcal{V}}\le C_pQ,\quad
 \|f_1-f_2\|_{\HH_1^1,p,\mathcal{V}}\le C_pQ\sigma,\quad
\big\|(\inf_{x\in[0,1]}f_2(x))^{-1}\big\|_{\R,p,\mathcal{V}}\le C_p,
\end{align}
where $Q\ge1$ is $\IF_0$-measurable and $0\le \sigma\le 1$ is deterministic. Then for every $\eta\in(0,1]$
\[
 \E\Big[
 \Big\|\mathscr R_{N,\eta}(f_1)
 -P_N\frac{\partial_x f_2}{f_2}\Big\|_H^2\mathbf1_{\mathcal{V}}
 \Big]
 \le C\big(1+\|Q\|_{4,\infty}^4\big)\left(\Psi(\sigma)+\eta^2 +\eta^{-2}\sigma^4\right).
\]
\end{lemma}

\begin{proof}
Define $Y:=    \mathbf1_{\mathcal{V}}\inf_{x\in[0,1]}f_2(x)+\mathbf1_{\mathcal{V}^c}$ and the events
\begin{align*}
\mathcal V_0&:=
\mathcal V\cap
\{\|f_1\|_{\mathcal H_1^1}=0\},\\
 \mathcal{V}_1&:=\mathcal{V}\cap
\{\|f_1\|_{\mathcal H_1^1}>0,~Y\ge2\eta\|f_1\|_{\HH^1_1},~ \|f_1-f_2\|_{\infty}\le Y/2\},\\
 \mathcal{V}_2&:=\mathcal{V}\cap
\{\|f_1\|_{\mathcal H_1^1}>0,~Y<2\eta\|f_1\|_{\HH^1_1}\},\\
 \mathcal{V}_3&:=\mathcal{V}\cap
\{\|f_1\|_{\mathcal H_1^1}>0,~Y\ge2\eta\|f_1\|_{\HH^1_1},~ \|f_1-f_2\|_{\infty}>Y/2\}.
\end{align*}
Thus $Y>0$ a.s., and we have
   \begin{align*}
   	\E\Big[\Big\|\mathscr R_{N,\eta}(f_1)-P_N\frac{\partial_x f_2}{f_2}\Big\|_H^2\mathbf1_{\mathcal{V}}\Big]&=  \E\Big[\Big\|P_N\frac{\partial_x f_2}{f_2}\Big\|_H^2\mathbf1_{\mathcal{V}_0}\Big]+\sum_{i=1}^3\E\Big[\Big\|P_{N}\frac{\partial_x f_1}{f_1\vee(\eta\|f_1\|_{\HH^1_1})}-P_N\frac{\partial_x f_2}{f_2}\Big\|_H^2\mathbf1_{\mathcal{V}_i}\Big]\\
   	&=:\mathbb{L}_0+\sum_{i=1}^3\mathbb{L}_i.
   \end{align*}
 
  \textbf{For the term $\mathbb{L}_0$}, since \(f_1=0\) in \(\mathcal H_1^1\) on \(\mathcal V_0\), and \(\mathscr R_{N,\eta}(0)=0\), the H\"older inequality and hypotheses in \eqref{conditions-f1f2} give
\begin{align*}
\mathbb L_0\le \E\big[(\inf_{x\in[0,1]}f_2(x))^{-2}\|f_2\|_{\HH^1_1}^2\mathbf1_{\mathcal{V}_0}\big]&\le \E\big[\big\|(\inf_{x\in[0,1]}f_2(x))^{-1}\big\|_{\R,4,\mathcal{V}}^{2}\|f_1-f_2\|_{\HH_1^1,4,\mathcal{V}}^{2}\big]\\
	 &\le C\sigma^2\E[Q^2]\le C\big(1+\|Q\|_{4,\infty}^4\big)\Psi(\sigma).
\end{align*}

\textbf{For the term $\mathbb{L}_1$}, since $Y\ge 2\eta\|f_1\|_{\HH^1_1}$ and $\|f_1-f_2\|_{\infty}\le Y/2$ on \(\mathcal V_1\), for every $x\in[0,1]$, 
   \begin{align*}
   	f_1(x)\ge f_2(x)-|f_1(x)-f_2(x)|\ge Y-\sup_{x\in[0,1]}|f_1(x)-f_2(x)|\ge Y-\|f_1-f_2\|_{\infty}\ge Y/2\ge\eta\|f_1\|_{\HH^1_1}.
   \end{align*}
 Thus we have
      \begin{align*}
    	\Big\|\frac{\partial_x f_1}{f_1}-\frac{\partial_x f_2}{f_2}\Big\|_H&=\Big\|\frac{\partial_x (f_1-f_2)f_2-\partial_x f_2(f_1-f_2)}{f_1f_2} \Big\|_H\\
    	&\le 2Y^{-1}\Big(\|f_1-f_2\|_{\HH^1_1}+\Big\|\frac{\partial_x f_2}{f_2}\Big\|_H\|f_1-f_2\|_{\infty}\Big).
    \end{align*} 
Using the embedding 
$\HH_1^1\hookrightarrow C([0,1];\R)$ we obtain 
    \begin{align*}
     \mathbb{L}_1=\E\Big[\Big\|P_N\Big(\frac{\partial_x  f_1}{f_1}-\frac{\partial_x f_2}{f_2}\Big)\Big\|_H^2\mathbf{1}_{\mathcal{V}_1}\Big]& \le 4\E\Big[ Y^{-2}\Big(\|f_1-f_2\|_{\HH^1_1}+\Big\|\frac{\partial_x f_2}{f_2}\Big\|_H\Big((C_E\|f_1-f_2\|_{\HH^1_1})\wedge \frac{Y}{2}\Big)\Big)^2\mathbf{1}_{\mathcal{V}_1}\Big].
   \end{align*}
Therefore, the H\"older inequality and hypotheses in \eqref{conditions-f1f2} yield  $\mathbb{L}_1\le C\E[Q^2(Q^2\sigma^2\wedge 1)].$ It is known that, by the layer cake representation, there exists a constant $C>0$ such that  
\begin{align}\label{ineq:tail1}
\E[Q^4\mathbf{1}_{\{Q\le R\}}]\le 4\int_0^{R}	y^{3}\mathbb P(Q>y)\,dy&\le 4\int_0^1 y^3 \,dy+4\big(\sup_{y\ge 1}y^4\mathbb P(Q>y)\big)\int_1^R y^{-1}\,dy\nonumber\\
&\le C(1+\|Q\|_{4,\infty}^4)(1+\ln R),
\end{align}
and similarly,
\begin{align}\label{ineq:tail2}
	\E[Q^2\mathbf{1}_{\{Q> R\}}]=R^2\mathbb P(Q> R)+2\int_R^\infty y\mathbb P(Q>y)\,dy&\le \big(\sup_{y\ge 1}y^4\mathbb P(Q>y)\big)\Big(R^{-2}+2\int_R^\infty y^{-3}dy\Big)\nonumber\\
	&\le 2(1+\|Q\|_{4,\infty}^4)R^{-2}.
\end{align}
Hence, for \(0<\sigma\le1\), splitting at $R=\sigma^{-1}$ and combining \eqref{ineq:tail1}--\eqref{ineq:tail2}, we obtain
 \begin{align*}
  	  \mathbb{L}_1\le C\sigma^2\E[Q^4\mathbf{1}_{\{Q\le \sigma^{-1}\}}]+C\E[Q^2\mathbf{1}_{\{Q> \sigma^{-1}\}}]\le C\big(1+\|Q\|_{4,\infty}^4\big)\sigma^2\left(1+0\vee\log(1/\sigma)\right)= C\bigl(1+\|Q\|_{4,\infty}^4\bigr)\Psi(\sigma).
  \end{align*}
  The case \(\sigma=0\) is immediate.

\textbf{For the term $\mathbb{L}_2$}, since $Y< 2\eta\|f_1\|_{\HH^1_1}$ on $\mathcal{V}_2$, we have 
\[
\mathcal V_2
\subset
\mathcal V\cap
\Big\{
\frac{\|f_1\|_{\mathcal H_1^1}}{Y}>(2\eta)^{-1}
\Big\}.
\]
It suffices to consider \(0<\eta\le1/2\) since the desired estimate is trivial for \(\eta\in(1/2,1]\). 
Using conditional Markov's inequality together with the hypotheses in \eqref{conditions-f1f2}, we obtain
\begin{align}\label{ineq:V2}
\mathbb P(\mathcal V_2)&\le \mathbb P\left(Q>(2\eta)^{-1}\right)+(2\eta)^8\E\Big[\E_0\big[
Y^{-8}\|f_1\|_{\mathcal H_1^1}^8
\mathbf1_{\mathcal V}
\big]
\mathbf1_{\{Q\le(2\eta)^{-1}\}}
\Big]\nonumber\\
&\le
\mathbb P\left(Q>(2\eta)^{-1}\right)
+
C(2\eta)^8
\E\left[
Q^8\mathbf1_{\{Q\le(2\eta)^{-1}\}}
\right]
\le C\bigl(1+\|Q\|_{4,\infty}^4\bigr)\eta^4.  
\end{align}
Here the last step follows from the weak-\(L^4\) tail bound and the corresponding truncated eighth-moment estimate similar to \eqref{ineq:tail1}.
By the contraction property of \(P_N\),
\begin{align*}
 \mathbb{L}_2&=\E\Big[\Big\|P_{N}\frac{\partial_x f_1}{f_1\vee(\eta\|f_1\|_{\HH^1_1})}-P_N\frac{\partial_x f_2}{f_2}\Big\|_H^2\mathbf{1}_{\mathcal{V}_2}\Big]\le 2\E\Big[\eta^{-2}\|f_1\|_{\HH^1_1}^{-2}\|\partial_x f_1\|_H^2\mathbf{1}_{\mathcal{V}_2}+\Big\|\frac{\partial_x f_2}{f_2}\Big\|_H^2\mathbf{1}_{\mathcal{V}_2}\Big].
\end{align*} 
Since
\(\|\partial_xf_1\|_H\le\|f_1\|_{\mathcal H_1^1}\) and \eqref{ineq:V2},
the first term satisfies
\[
\E\big[
\eta^{-2}\|f_1\|_{\mathcal H_1^1}^{-2}
\|\partial_xf_1\|_H^2\mathbf1_{\mathcal V_2}
\big]
\le
\eta^{-2}\mathbb P(\mathcal V_2)
\le
C\bigl(1+\|Q\|_{4,\infty}^4\bigr)\eta^2.
\]
For the second term, we split according to
\(Q\le(2\eta)^{-1}\) and \(Q>(2\eta)^{-1}\). By conditional H\"older's inequality,
\[
\begin{aligned}
&\E\Big[
\Big\|\frac{\partial_xf_2}{f_2}\Big\|_H^2
\mathbf1_{\mathcal V_2}
\mathbf1_{\{Q\le(2\eta)^{-1}\}}
\Big]\le
\E\Big[
\Big(
\E_0\Big[
\Big\|\frac{\partial_xf_2}{f_2}\Big\|_H^6
\mathbf1_{\mathcal V}
\Big]
\Big)^{1/3}
\big(
\E_0[\mathbf1_{\mathcal V_2}]
\big)^{2/3}
\mathbf1_{\{Q\le(2\eta)^{-1}\}}
\Big].
\end{aligned}
\]
Conditional Markov's inequality and hypotheses in \eqref{conditions-f1f2} additionally give
\[
\E_0[\mathbf1_{\mathcal V_2}]
\le
(2\eta)^6
\E_0\left[
Y^{-6}\|f_1\|_{\mathcal H_1^1}^6
\mathbf1_{\mathcal V}
\right]
\le
C(2\eta)^6Q^6.
\]
Consequently, we obtain
\[
\begin{aligned}
&\E\Big[
\Big\|\frac{\partial_xf_2}{f_2}\Big\|_H^2
\mathbf1_{\mathcal V_2}
\mathbf1_{\{Q\le(2\eta)^{-1}\}}
\Big]\le
C\eta^4
\E\left[
Q^6\mathbf1_{\{Q\le(2\eta)^{-1}\}}
\right]
\le
C\bigl(1+\|Q\|_{4,\infty}^4\bigr)\eta^2.
\end{aligned}
\]
On the other hand, by hypotheses in \eqref{conditions-f1f2} and \eqref{ineq:tail2}, it holds that
\[
\begin{aligned}
&\E\Big[
\Big\|\frac{\partial_xf_2}{f_2}\Big\|_H^2
\mathbf1_{\mathcal V}
\mathbf1_{\{Q>(2\eta)^{-1}\}}
\Big]\le
C\E\left[
Q^2\mathbf1_{\{Q>(2\eta)^{-1}\}}
\right]
\le
C\bigl(1+\|Q\|_{4,\infty}^4\bigr)\eta^2.
\end{aligned}
\]
Combining the above estimates yields $\mathbb L_2
\le
C\bigl(1+\|Q\|_{4,\infty}^4\bigr)\eta^2.$

\textbf{For the term $\mathbb{L}_3$},  since $\|f_1-f_2\|_{\infty}>Y/2$ on $\mathcal{V}_3$, we have $\mathcal V_3
\subset
\mathcal V\cap
\big\{
\frac{\|f_1-f_2\|_\infty}{\sigma Y}>(2\sigma)^{-1}
\big\}.$  
Conditional Markov's inequality gives, for $0<\sigma\le 1/2$, 
\[
\begin{aligned}
\mathbb P\Big(
\mathcal V\cap
\Big\{
\frac{\|f_1-f_2\|_\infty}{\sigma Y}>(2\sigma)^{-1}
\Big\}
\Big)
&\le
\mathbb P(Q>(2\sigma)^{-1})
+
2^8
\E\left[
\E_0\left[
\|f_1-f_2\|_\infty^8Y^{-8}
\mathbf1_{\mathcal V}
\right]
\mathbf1_{\{Q\le (2\sigma)^{-1}\}}
\right].
\end{aligned}
\]
By the embedding  
\(\mathcal H_1^1\hookrightarrow C([0,1];\mathbb R)\), conditional H\"older's inequality, and hypotheses in \eqref{conditions-f1f2},
\[
\begin{aligned}
\E_0\big[
\|f_1-f_2\|_\infty^8Y^{-8}\mathbf1_{\mathcal V}
\big]\le C
\left(\E_0\big[
\|f_1-f_2\|_{\mathcal H_1^1}^{16}
\mathbf1_{\mathcal V}
\big]
\right)^{1/2}
\left(
\E_0\left[
Y^{-16}\mathbf1_{\mathcal V}
\right]
\right)^{1/2}\le
C Q^8\sigma^8.
\end{aligned}
\]
Hence it follows that 
\begin{align}\label{proba-v3}
\mathbb P\Big(
\mathcal V\cap
\Big\{
\frac{\|f_1-f_2\|_\infty}{\sigma Y}>(2\sigma)^{-1}
\Big\}
\Big)
&\le
\mathbb P\bigl(Q>(2\sigma)^{-1}\bigr)+C\sigma^8
\E\left[
Q^8\mathbf1_{\{Q\le(2\sigma)^{-1}\}}
\right]\le
C\bigl(1+\|Q\|_{4,\infty}^4\bigr)\sigma^4,
\end{align}
where the last step follows from the weak-\(L^4\) tail bound and the corresponding truncated eighth-moment estimate. This implies that $\mathbb P(\mathcal V_3)
\le
C\bigl(1+\|Q\|_{4,\infty}^4\bigr)\sigma^4.$ 
Similar to estimates for term $\mathbb L_2$, 
\begin{align*}
\mathbb L_3
&=
\E\Big[
\Big\|
P_N\frac{\partial_xf_1}
{f_1\vee(\eta\|f_1\|_{\mathcal H_1^1})}
-
P_N\frac{\partial_xf_2}{f_2}
\Big\|_H^2
\mathbf1_{\mathcal V_3}
\Big]\le
C\eta^{-2}\mathbb P(\mathcal V_3)
+
C\E\Big[
\Big\|\frac{\partial_xf_2}{f_2}\Big\|_H^2
\mathbf1_{\mathcal V_3}
\Big].
\end{align*}
It remains to estimate the second term. We split this term into
\begin{align*}
\E\Big[
\Big\|\frac{\partial_xf_2}{f_2}\Big\|_H^2
\mathbf1_{\mathcal V_3}
\Big]&=
\E\Big[
\Big\|\frac{\partial_xf_2}{f_2}\Big\|_H^2
\mathbf1_{\mathcal V_3}
\mathbf1_{\{Q\le (2\sigma)^{-1}\}}
\Big]
+
\E\Big[
\Big\|\frac{\partial_xf_2}{f_2}\Big\|_H^2
\mathbf1_{\mathcal V_3}
\mathbf1_{\{Q>(2\sigma)^{-1}\}}
\Big].
\end{align*}
Since $\mathcal V_3
\subset
\mathcal V\cap
\{\|f_1-f_2\|_\infty>Y/2\},$ conditional Markov's inequality, the embedding 
$\mathcal H_1^1\hookrightarrow C([0,1];\mathbb R)$,
and hypotheses in \eqref{conditions-f1f2}, imply
\begin{align*}
\E_0[\mathbf1_{\mathcal V_3}]
&\le
2^4\E_0\Big[
Y^{-4}\|f_1-f_2\|_\infty^4
\mathbf1_{\mathcal V}
\Big]\le
C
\big(
\E_0\big[
Y^{-8}\mathbf1_{\mathcal V}
\big]
\big)^{1/2}
\big(
\E_0\big[
\|f_1-f_2\|_{\mathcal H_1^1}^8
\mathbf1_{\mathcal V}
\big]
\big)^{1/2} \le
CQ^4\sigma^4.
\end{align*}
Hence, conditional H\"older's inequality gives
\begin{align*}
\E\Big[
\Big\|\frac{\partial_xf_2}{f_2}\Big\|_H^2
\mathbf1_{\mathcal V_3}
\mathbf1_{\{Q\le (2\sigma)^{-1}\}}
\Big]
&\le
\E\Big[
\Big(
\E_0\Big[
\Big\|\frac{\partial_xf_2}{f_2}\Big\|_H^8
\mathbf1_{\mathcal V}
\Big]
\Big)^{1/4}
\Big(
\E_0[\mathbf1_{\mathcal V_3}]
\Big)^{3/4}
\mathbf1_{\{Q\le (2\sigma)^{-1}\}}
\Big]\le
C\sigma^3
\E[
Q^5\mathbf1_{\{Q\le (2\sigma)^{-1}\}}].
\end{align*}
Since $\E[
Q^5\mathbf1_{\{Q\le (2\sigma)^{-1}\}}
]
\le
C(1+\|Q\|_{4,\infty}^4)(2\sigma)^{-1},$ and therefore 
\[
\E\Big[
\Big\|\frac{\partial_xf_2}{f_2}\Big\|_H^2
\mathbf1_{\mathcal V_3}
\mathbf1_{\{Q\le (2\sigma)^{-1}\}}
\Big]
\le
C\bigl(1+\|Q\|_{4,\infty}^4\bigr)\sigma^2.
\]
For the remaining part, by hypotheses in \eqref{conditions-f1f2} and~\eqref{ineq:tail2},
\begin{align*}
\E\Big[
\left\|\frac{\partial_xf_2}{f_2}\right\|_H^2
\mathbf1_{\mathcal V_3}
\mathbf1_{\{Q>(2\sigma)^{-1}\}}
\Big]&\le
\E\Big[
\E_0\Big[
\Big\|\frac{\partial_xf_2}{f_2}\Big\|_H^2
\mathbf1_{\mathcal V}
\Big]
\mathbf1_{\{Q>(2\sigma)^{-1}\}}
\Big]\le
C\E\Big[
Q^2\mathbf1_{\{Q>(2\sigma)^{-1}\}}
\Big]\le
C\bigl(1+\|Q\|_{4,\infty}^4\bigr)\sigma^2.
\end{align*}
Combining the above estimates yields
$\mathbb L_3
\le
C\bigl(1+\|Q\|_{4,\infty}^4\bigr)
\big(\eta^{-2}\sigma^4+\sigma^2\big).$ The case \(\sigma=0\) is immediate.

Consequently, the estimates for terms $\mathbb{L}_0$--$\mathbb{L}_3$ yield 
\[
 \E\Big[
 \Big\|\mathscr R_{N,\eta}(f_1)
 -P_N\frac{\partial_x f_2}{f_2}\Big\|_H^2\mathbf1_{\mathcal{V}}
 \Big]
 \le C\big(1+\|Q\|_{4,\infty}^4\big)\left(\Psi(\sigma)+\eta^2 +\eta^{-2}\sigma^4\right).
\]
This completes the proof.
\end{proof}

\subsection{Proof of Theorem~\ref{thm:sharp-main}}

\begin{proof}

Fix $0<\varsigma<\frac12-2\alpha$ and recall that $\delta=1\wedge\bigl(N^{-\varsigma}+\sqrt{\tau}\bigr).$  It suffices to consider the nontrivial regime $N^{-2}+N^{-\varsigma}+\sqrt{\tau}\le 1/2$ so that $\delta=N^{-\varsigma}+\sqrt{\tau}\le 1/2$. Indeed, if \(N^{-2}+N^{-\varsigma}+\sqrt{\tau}>1/2\), then $N^{-\varsigma}+\sqrt{\tau}>1/4$, and then the desired estimate follows directly from the uniform second-moment bounds for solution \(u\) and numerical approximation \(u_m^{N,\delta}\).

Set $\ell_N:=\frac{2\log(eN)}{\pi^2N^2}.$ We first consider the initial time $m=0$. 
The definition of \(\mathcal J_\varsigma\) gives
\[
\E\bigl[
\|u(0)-u_0^{N,\delta}\|_H^2
\bigr]
=
\E\bigl[
\|u_0-\mathcal R_{N,N^{-\varsigma}}(U^N_0)\|_H^2
\bigr]
\le
\mathcal J_\varsigma N^{-2\varsigma}.
\]

It therefore remains to estimate the error at the positive grid points
\(t_m\), \(m=1,\ldots,M\). The subsequent analysis is divided into two cases according to whether \(t_m\) lies beyond or within the initial layer:
\[
t_m\ge\ell_N
\qquad\text{and}\qquad
t_m<\ell_N.
\]



\subsubsection{Case 1. $t_m\ge \ell_N$}

Set $\mathcal A_N:=\{\|U_0\|_{\mathcal H_1^1}\le N\}.$ We decompose
\[
\E[\|u(t_m)-u_m^{N,\delta}\|_H^2]
=
\E\!\left[\|u(t_m)-u_m^{N,\delta}\|_H^2\mathbf1_{\mathcal A_N}\right]
+
\E\!\left[\|u(t_m)-u_m^{N,\delta}\|_H^2\mathbf1_{\mathcal A_N^c}\right].
\]

We first estimate the error on \(\mathcal A_N^c\). Using the basic inequality $\|a-b\|_H^2\le 2\|a\|_H^2+2\|b\|_H^2$, together with  \eqref{eq:sharp-u-moments} and
\(\|\mathscr R_{N,\delta}(f)\|_H\le\delta^{-1}\), we obtain
\begin{align}\label{ineq:othercase}
\E\!\left[
\|u(t_m)-u_m^{N,\delta}\|_H^2\mathbf1_{\mathcal A_N^c}\right]&\le 2\E[\|u(t_m)\|_{L_0^2(H)}^2\mathbf1_{\mathcal A_N^c}]+2\delta^{-2}\mathbb P(\mathcal A_N^c)\nonumber\\
&\le C\E\big[
\|U_0\|_{\mathcal H_1^1}^2
\mathbf1_{\mathcal A_N^c}\big]+
C\delta^{-2}N^{-2}
\E\big[
\|U_0\|_{\mathcal H_1^1}^2
\mathbf1_{\mathcal A_N^c}
\big]\nonumber\\
&\le
C(1+\|\,\|U_0\|_{\mathcal H_1^1}\|_{4,\infty}^4)
(N^{-2}+\delta^{-2}N^{-4})\le
C(1+\|\,\|U_0\|_{\mathcal H_1^1}\|_{4,\infty}^4)
N^{-2\varsigma},
\end{align}
where we used \eqref{ineq:tail2} and $\delta^{-2}N^{-4}
\le N^{2\varsigma-4}\le N^{-2\varsigma}$ for $\varsigma\le1$.

To estimate the error on $\mathcal A_N$,  we then insert the spectral projection \(P_Nu(t_m)\) and decompose
\begin{align*}
 \E\big[ \|u(t_m)-u^{N,\delta}_m\|_H^2\mathbf{1}_{\mathcal A_N}\big]&\le 2\big(\E\big[ \|u(t_m)-P_Nu(t_m)\|_H^2\mathbf{1}_{\mathcal A_N}\big]+\E\big[ \|P_Nu(t_m)-u^{N,\delta}_m\|_H^2\mathbf{1}_{\mathcal A_N}\big]\big)=:2(\mathbb{T}_1+\mathbb{T}_2). 
\end{align*}

\textbf{For the term $\mathbb{T}_1$}, the mild formulation gives 
\begin{align*}
\E\big[
\|u(t_m)-P_Nu(t_m)\|_H^2
\mathbf1_{\mathcal A_N}
\big]&\le
3\E[\|(\mathrm I-P_N)S(t_m)u_0\|_H^2]+
3\E\Big[\Big\|
\sum_{k\ge N+1}\lambda_k^\alpha
\int_0^{t_m}S(t_m-s)e_k\,d\beta_k(s)
\Big\|_H^2\Big]\\
&\qquad+
3\Big(
\int_0^{t_m}
\big(\E\big[\big\|
(\mathrm I-P_N)S(t_m-s)\partial_x(u(s)^2)
\mathbf1_{\mathcal A_N}
\big\|_H^2\big]\big)^{1/2}
\,ds\Big)^2. 
\end{align*}
For the initial value, since
\(\mathcal R_{N,N^{-\varsigma}}(U^N_0)\in H_N\), the best-approximation property of \(P_N\) gives
\begin{align*}
  \|(\mathrm I-P_N)S(t_m)u_0\|_{L^2(\Omega;H)}^2
\le C_\varsigma N^{-2\varsigma}
\big(
N^{2\varsigma}\|(\mathrm I-P_N)u_0\|_{L^2(\Omega;H)}^2
\big)\le C_\varsigma\mathcal{J}_{\varsigma}N^{-2\varsigma}. 
\end{align*}
For the stochastic convolution, It\^o's isometry and \(S(t)e_k=e^{-\lambda_kt}e_k\) yield
\begin{align*}
  \E\Big[\Big\|
\sum_{k\ge N+1}\lambda_k^\alpha
\int_0^{t_m}S(t_m-s)e_k\,d\beta_k(s)
\Big\|_H^2\Big]=
\frac12\sum_{k\ge N+1}
\lambda_k^{2\alpha-1}
\bigl(1-e^{-2\lambda_kt_m}\bigr)
\le C_\alpha N^{-1+4\alpha}.  
\end{align*}
Lemma~\ref{lem:heat-tail} yields, for any $g\in H$ and all $t\in(0,T]$, $\|(\mathrm I-P_N)S(t)\partial_x(g^2)\|_H
 \le C_Tt^{-3/4}e^{-cN^2t}\|g\|_H^2.$ Thus, for the nonlinear term, we obtain
\[
\begin{aligned}
&~~\quad \Big(
\int_0^{t_m}
\big(\E\big[\big\|
(\mathrm I-P_N)S(t_m-s)\partial_x(u(s)^2)
\mathbf1_{\mathcal A_N}
\big\|_H^2\big]\big)^{1/2}\,ds
\Big)^2\\
&\le C_T\Big(\int_0^{t_m}(t_m-s)^{-3/4}e^{-cN^2(t_m-s)}\Bigl(
\E[\|u(s)\|_H^4\mathbf1_{\mathcal A_N}]
\Bigr)^{1/2}\,ds\Big)^2\\
&\le
C_TN^{-1}
\sup_{0\le s\le t_m}
\E[\|u(s)\|_H^4\mathbf1_{\mathcal A_N}].
\end{aligned}
\]
Results \eqref{eq:sharp-u-moments} and \eqref{ineq:tail1}, and the fact $N^{-1}(1+\ln N)\le C_{\varsigma}N^{-2\varsigma}$ for $0<\varsigma<\frac12$  give
\[
\begin{aligned}
N^{-1}
\sup_{0\le s\le t_m}
\E[\|u(s)\|_H^4\mathbf1_{\mathcal A_N}]\le CN^{-1}\bigl(
1+\|\|U_0\|_{\mathcal H_1^1}\|_{4,\infty}^4
\bigr)(1+\ln N)\le
C_{\varsigma,T}
\bigl(
1+\|\|U_0\|_{\mathcal H_1^1}\|_{4,\infty}^4
\bigr)N^{-2\varsigma}.
\end{aligned}
\]

Combining the estimates for the three subterms gives, for $\varsigma<1/2-2\alpha$,
\begin{align}\label{ineq:T1} 
\mathbb T_1&\le
C_\alpha N^{-1+4\alpha}+C_{\varsigma,T}
\bigl(1+\mathcal J_\varsigma+\big\|\|U_0\|_{\mathcal{H}_1^1}\big\|_{4,\infty}^4\bigr)N^{-2\varsigma}\le C_{\alpha,\varsigma,T}\big(1+\mathcal J_\varsigma+\big\|\|U_0\|_{\mathcal{H}_1^1}\big\|_{4,\infty}^4\big)N^{-2\varsigma}.
\end{align}

\textbf{For the term $\mathbb{T}_2$}, since
$t_m\ge\ell_N$, we have
\[
 e^{-\lambda_{N+1}t_m}
 \le \exp\!\Big[-2\frac{(N+1)^2}{N^2}\log(eN)\Big]
 \le (eN)^{-2}\le N^{-2}.
\]
Applying \eqref{eq:sharp-linear-error} with $\gamma=1$, for $\varsigma<a<\frac12-2\alpha$, we obtain that for all $m\in\{1,2,\ldots,M\}$ 
  \begin{align*}
  \|U(t_m)-U^N_m\|_{L^p_0(\HH^{1}_1)}\le C_{p}\|U_0\|_{\HH^{1}_1}(N^{-2}+N^{-a}+\sqrt{\tau}). 
\end{align*}  
By \eqref{eq:sharp-moments} and \eqref{eq:sharp-inverse-moments}, the hypotheses on Lemma~\ref{lem:error-f1f2} are satisfied with $\mathcal{V}=\Omega$, $Q=\|U_0\|_{\HH^1_1}$, $\eta=\delta$, and $\sigma=N^{-2}+N^{-a}+\sqrt{\tau}$. If we take $f_1=U^N_m$ and $f_2=U(t_m)$ in Lemma~\ref{lem:error-f1f2}, by Lemma~\ref{lem:sharp-existence}, it holds that  
\begin{align*}
   \mathbb T_2&=\E\big[ \|P_Nu(t_m)-u^{N,\delta}_m\|_H^2\mathbf{1}_{\mathcal{A}_N}\big]\le \E\Big[\Big\|P_N\frac{\partial_x U(t_m)}{U(t_m)}-\mathscr{R}_{N,\delta}(U^N_m)\Big\|_H^2\Big]\nonumber\\
   &\le C\big(1+\big\|\|U_0\|_{\mathcal{H}_1^1}\big\|_{4,\infty}^4\big)\left(\Psi(N^{-2}+N^{-a}+\sqrt{\tau})+\delta^2 +\delta^{-2}(N^{-2}+N^{-a}+\sqrt{\tau})^4\right).
\end{align*}
Recall that $\delta=N^{-\varsigma}+\sqrt{\tau}$ and $\varsigma<a<1/2-2\alpha.$ Since \(N^{-2}+N^{-a}+\sqrt{\tau}\le C\delta\), we have 
\[
\delta^2+\delta^{-2}
\bigl(N^{-2}+N^{-a}+\sqrt{\tau}\bigr)^4
\le
C\bigl(N^{-2\varsigma}+\tau\bigr).
\]
Moreover, by the definition of \(\Psi\), we obtain that $\Psi(N^{-2}+N^{-a}+\sqrt{\tau})
\le
C_{a,\varsigma}\bigl(
N^{-2\varsigma}
+
\tau\bigl(1+0\vee\ln(1/\tau)\bigr)
\bigr).$ Consequently,
\begin{align}\label{ineq:T2} 
   \mathbb T_2\le C_{a,\varsigma}\big(1+\big\|\|U_0\|_{\mathcal{H}_1^1}\big\|_{4,\infty}^4\big)(N^{-2\varsigma}+\tau(1+0\vee\log (1/\tau))).
\end{align}

Consequently, combining the results \eqref{ineq:othercase}, \eqref{ineq:T1}, and \eqref{ineq:T2}, there exists a constant $C_{\alpha,\varsigma,T}>0$ such that
 \begin{align*}
   \E[\|u(t_m)-u^{N,\delta}_m\|_H^2]\le C_{\alpha,\varsigma,T}\big(1+\mathcal{J}_{\varsigma}+\big\|\|U_0\|_{\mathcal{H}_1^1}\big\|_{4,\infty}^4\big)\big(N^{-2\varsigma}+\tau(1+0\vee\log (1/\tau))\big). 
 \end{align*}
We complete the proof of Theorem~\ref{thm:sharp-main} in the case that $t_m\ge\ell_N$.

\subsubsection{Case 2. $t_m<\ell_N$}

Let \(m\in\{1,2,\ldots,M\}\) satisfy $t_m<\ell_N<1$. Since $U_0\ge1$, 
$\|u_0\|_H=\|\partial_x U_0/U_0\|_H\le \|U_0\|_{\HH_1^1}$.
For all $x\in[0,1]$, by the Cauchy--Schwarz inequality, there exists a constant $C_1\in(0,\infty)$ such that
\begin{align*}
	|(\mathrm I-\Pi_N)U_0(x)|&\le \sum_{k\ge N+1}|\langle U_0, h_k\rangle_H| |h_k(x)|\le \sqrt{2}\sum_{k\ge N+1}k|\langle U_0, h_k\rangle_H| k^{-1}\\
	&\le\sqrt{2}\Big(\sum_{k\ge N+1}k^2|\langle U_0, h_k\rangle_H|^2\Big)^{1/2}\Big(\sum_{k\ge N+1}k^{-2}\Big)^{1/2}\le C_1N^{-1/2}\|U_0\|_{\HH^1_1}.
\end{align*}
Choose $c_0>0$ such that $C_1c_0\le1/2$ and define
\[
\mathcal{G}_{N,c_0}:=\{\|U_0\|_{\HH_1^1}\le c_0\sqrt{N}\}
\]
Then, on \(\mathcal G_{N,c_0}\), 
\begin{align*}
 U^N_0(x)\ge U_0(x)-|(\mathrm I-\Pi_N)U_0(x)|\ge U_0(x)-C_1N^{-1/2}c_0\sqrt{N}\ge 1-C_1c_0\ge 1/2.
\end{align*}
As detailed in \eqref{Neumann-semigroup}, for all $t\in[0,T]$ and all $x\in[0,1]$, it holds that
\begin{align}\label{ineq:UN0}
	(S_1(t)U^N_0)(x)=\E_R[U^N_0(R^x(t))]\ge 1/2\qquad \mbox{on}\quad \mathcal G_{N,c_0}.
\end{align}
In particular, all quotients involving \(S_1(t)U_0^N\) below are well defined on this event.

Now we decompose
\begin{align*}
   \E\big[\|u(t_m)-u^{N,\delta}_m\|_H^2\big]=\E\big[\|u(t_m)-u^{N,\delta}_m\|_H^2\mathbf{1}_{\mathcal{G}_{N,c_0}}\big] +\E\big[\|u(t_m)-u^{N,\delta}_m\|_H^2\mathbf{1}_{\mathcal{G}_{N,c_0}^c}\big]=:\mathbb{S}_1+\mathbb{S}_2. 
\end{align*}

\textbf{For the term $\mathbb{S}_1$},  by Lemma~\ref{lem:sharp-existence}, using \(u=\partial_xU/U\), we further split \(\mathbb S_1\) as 
\begin{align*}
\mathbb{S}_1\le 3\bigl(\mathbb{S}_{11}+\mathbb{S}_{12}+\mathbb{S}_{13}\bigr),
\end{align*}
where
\begin{align*}
  &\mathbb{S}_{11}:= \E\Big[\Big\|\frac{\partial_x U(t_m)}{U(t_m)}-\frac{\partial_x (S_1(t_m)U_0)}{S_1(t_m)U_0}\Big\|_H^2\mathbf{1}_{\mathcal{G}_{N,c_0}}\Big],\\
  &\mathbb{S}_{12}:=\E\Big[\Big\|\frac{\partial_x (S_1(t_m)U_0)}{S_1(t_m)U_0}-P_N\frac{\partial_x (S_1(t_m)U^N_0)}{S_1(t_m)U^N_0}\Big\|_H^2\mathbf{1}_{\mathcal{G}_{N,c_0}}\Big],\\
  &\mathbb{S}_{13}:=\E\Big[\Big\|P_N\frac{\partial_x (S_1(t_m)U^N_0)}{S_1(t_m)U^N_0}-\mathscr{R}_{N,\delta}(U^{N}_m)\Big\|_H^2\mathbf{1}_{\mathcal{G}_{N,c_0}}\Big]. 
\end{align*}

\textit{1. Estimate of the subterm $\mathbb{S}_{11}$.} Since $U_0\ge 1$, for all $t\in[0,T]$ and all $x\in[0,1]$, we have 
\begin{align}\label{ineq:U0}
	(S_1(t)U_0)(x)=\E_R[U_0(R^x(t))]\ge 1.
\end{align}
Using \eqref{eq:sharp-inverse-moments} and \eqref{eq:sharp-start-exact} with $\gamma=1$ yields, for all $t\in(0,1\wedge T]$ and $\varsigma<a<\frac12-2\alpha$,
\begin{align*}
  &\Big\|\frac{\partial_x U(t)}{U(t)}-\frac{\partial_x (S_1(t)U_0)}{S_1(t)U_0}\Big\|_{L^2_0(H)}^2=\Big\|\frac{\partial_x (U(t)-S_1(t)U_0)}{U(t)}+\frac{\partial_x (S_1(t)U_0)(S_1(t)U_0-U(t))}{U(t)(S_1(t)U_0)}\Big\|_{L^2_0(H)}^2\\
  \le& C\E_0\Big[\bigl(\inf_{x\in[0,1]}U(t,x)\bigr)^{-2}\|\partial_x(U(t)-S_1(t)U_0)\|_{H}^2\Big]\\
  &+C\E_0\Big[\bigl(\inf_{x\in[0,1]}U(t,x)\bigr)^{-2}\|\partial_x  (S_1(t)U_0)\|_H^2\|U(t)-S_1(t)U_0\|_{C([0,1];\R)}^2\Big]\\
  \le& C\|U_0\|_{\HH^1_1}^2 \E_0\Big[\bigl(\inf_{x\in[0,1]}U(t,x)\bigr)^{-2}\|U(t)-S_1(t)U_0\|_{\HH^1_1}^2\Big]\le C\|U_0\|_{\HH^1_1}^{4}t^{a}.
\end{align*}
Hence, for the $\mathcal F_0$-measurable event \(\mathcal G_{N,c_0}\), 
the tower property and 
\eqref{ineq:tail1} yield
\begin{align*}
   \mathbb{S}_{11}
   &=\E\Big[\Big\|u(t_m)-\frac{\partial_x (S_1(t_m)U_0)}{S_1(t_m)U_0}\Big\|_{L^2_0(H)}^2\mathbf{1}_{\mathcal G_{N,c_0}}\Big]\le Ct_m^{a}\E[\|U_0\|_{\HH^1_1}^4\mathbf{1}_{\mathcal G_{N,c_0}}]\le C\ell_N^{a}(1+\big\|\|U_0\|_{\mathcal{H}_1^1}\big\|_{4,\infty}^4)(1+\ln N).
\end{align*}
Recalling that $\ell_N=\frac{2\log(eN)}{\pi^2N^2}$ and $\varsigma<a$, we have $\ell_N^a(1+\ln N)
\le C_{a,\varsigma}N^{-2\varsigma}.$ Consequently,
\[
\mathbb S_{11}
\le
C_{a,\varsigma}
\bigl(1+\big\|\|U_0\|_{\mathcal{H}_1^1}\big\|_{4,\infty}^4\bigr)
N^{-2\varsigma}.
\]

\textit{2. Estimate of the subterm $\mathbb{S}_{12}$.} For $t\in[0,T]$, set 
\[
H_t:=\frac{\partial_x (S_1(t)U_0)}{S_1(t)U_0}\qquad \mbox{and}\qquad \mathbb{X}_t^N:=\frac{\partial_x (S_1(t)U^N_0)}{S_1(t)U^N_0}.
\]
 Their initial values are \(u_0\) and \(\hat u_0^N:=\partial_xU_0^N/U^N_0\), respectively. Applying Lemma~\ref{lem:deterministic-HC} to \(H_t\) and \(\mathbb X_t^N\) yields
\[
 \|H_t-P_N\mathbb{X}_t^N\|_H
 \le \|u_0-P_N\hat{u}^N_0\|_H
 +Ct^{1/4}\bigl(\|u_0\|_H^2+\|\hat{u}^N_0\|_H^2\bigr).
\]
Since $U_0\ge 1$ and $U^N_0\ge 1/2$, on $\mathcal G_{N,c_0}$, it holds that $\|u_0\|_H\le \|U_0\|_{\HH_1^1}$ and $\|\hat{u}^N_0\|_H\le2\|U_0\|_{\HH_1^1}$. Define $\mathcal{B}_N:=\{N^{-\varsigma}\|U_0\|_{\mathcal H_1^1}\le 1/2\}$. On $\mathcal G_{N,c_0}\cap\mathcal{B}_N$, we have $N^{-\varsigma}
\le\frac12$, and hence 
\[
U^N_0\ge
N^{-\varsigma},\qquad  P_N\hat{u}^N_0=\mathcal R_{N,N^{-\varsigma}}(U_0^N).
\]
On $\mathcal G_{N,c_0}
\cap
\mathcal{B}_N^c,$ the bounds $S_1(t)U_0\ge 1$ and $S_1(t)U_0^N\ge 1/2$ give
\[
\|H_t-P_N\mathbb{X}_t^N\|_H^2\le (\|H_t\|_{H}+\|\mathbb{X}_t^N\|_H)^2\le 9\|U_0\|_{\HH^1_1}^2.
\]
Therefore by \eqref{ineq:tail1} and \eqref{ineq:tail2}, and the fact that $\ell_N^{1/2}(1+\ln N)
\le C_{\varsigma}N^{-2\varsigma}$ for every \(\varsigma<1/2\), we obtain
\begin{align*}
\mathbb{S}_{12}
 &=\E\big[\|H_{t_m}-P_N\mathbb{X}_{t_m}^N\|_H^2\mathbf{1}_{\mathcal{G}_{N,c_0}\cap \mathcal{B}_N^c}\big]+\E\big[\|H_{t_m}-P_N\mathbb{X}_{t_m}^N\|_H^2\mathbf{1}_{\mathcal{G}_{N,c_0}\cap \mathcal{B}_N}\big] \\
 &\le C\E[\|U_0\|_{\HH^1_1}^2\mathbf{1}_{\mathcal B_N^c}]+Ct_m^{1/2}\E\big[\|U_0\|_{\HH^1_1}^4\mathbf{1}_{\mathcal{G}_{N,c_0}}\big]+C\E[\|u_0-\mathcal R_{N,N^{-\varsigma}}(U_0^N)\|_H^2\mathbf{1}_{\mathcal{G}_{N,c_0}\cap \mathcal{B}_N}]
 \\
 &\le C\big(1+\big\|\|U_0\|_{\mathcal{H}_1^1}\big\|_{4,\infty}^4\big)N^{-2\varsigma}+C\ell_N^{1/2}\big(1+\big\|\|U_0\|_{\mathcal{H}_1^1}\big\|_{4,\infty}^4\big)(1+\ln N)+CN^{-2\varsigma}\mathcal{J}_{\varsigma}\\
 &\le C_{a,\varsigma}\big(1+\mathcal{J}_{\varsigma}+\big\|\|U_0\|_{\mathcal{H}_1^1}\big\|_{4,\infty}^4\big)N^{-2\varsigma}.
\end{align*}

\textit{3. Estimate of the subterm $\mathbb{S}_{13}$.} 
By \eqref{eq:sharp-start-grid} and \eqref{ineq:UN0}, the hypotheses on Lemma~\ref{lem:error-f1f2} are satisfied with $\mathcal{V}=\mathcal{G}_{N,c_0}$,  $f_1=U^{N}_m$, $f_2=S_1(t_m)U^N_0$, $Q=\|U_0\|_{\HH^1_1}$, $\eta=\delta$,  and $\sigma=t_m^{a/2}$, $a\in(0,1/2-2\alpha)$. This implies that
\begin{align*}
   \mathbb{S}_{13}&= \E\Big[\Big\|P_N\frac{\partial_x (S_1(t_m)U^N_0)}{S_1(t_m)U^N_0}-\mathscr{R}_{N,\delta}(U^{N}_m)\Big\|_H^2 \mathbf{1}_{\mathcal{G}_{N,c_0}}\Big]\le C\big(1+\big\|\|U_0\|_{\mathcal{H}_1^1}\big\|_{4,\infty}^4\big)\left(\Psi(t_m^{a/2})+\delta^2 +\delta^{-2}t_m^{2a}\right).
\end{align*}
Since \(t_m<\ell_N<1\) and \(\delta=N^{-\varsigma}+\sqrt{\tau}\), we have $ \delta^{-2}t_m^{2a}\le \delta^{-2}\ell_N^{2a}
\le C_{a,\varsigma}N^{-2\varsigma}$ and $\Psi(t_m^{a/2})
\le
C_a\ell_N^a(1+|\ln\ell_N|)\le C_{a,\varsigma}N^{-2\varsigma}$. Hence for $\varsigma<a<1/2-2\alpha$
\begin{align*}
   \mathbb{S}_{13}
   &\le C_{a,\varsigma}\big(1+\big\|\|U_0\|_{\mathcal{H}_1^1}\big\|_{4,\infty}^4\big)(N^{-2\varsigma}+\tau).
\end{align*}

Combining the estimates of the subterms $\mathbb{S}_{11}$, $\mathbb{S}_{12}$, and $\mathbb{S}_{13}$ we obtain 
\begin{align*}
   \mathbb{S}_1\le  C_{a,\varsigma}\big(1+\mathcal{J}_{\varsigma}+\big\|\|U_0\|_{\mathcal{H}_1^1}\big\|_{4,\infty}^4\big)(N^{-2\varsigma}+\tau).
\end{align*}

\textbf{For the term $\mathbb{S}_2$}, recall that $\mathcal G_{N,c_0}^c
=
\{\|U_0\|_{\mathcal H_1^1}>c_0\sqrt N\}.$ Using \(\|a-b\|_H^2\le2\|a\|_H^2+2\|b\|_H^2\), together with the fact \(\|\mathscr R_{N,\delta}(f)\|_H\le\delta^{-1}\), we obtain
\begin{align*}
   \mathbb{S}_2=\E\big[\|u(t_m)-u^{N,\delta}_m\|_H^2\mathbf{1}_{\mathcal{G}_{N,c_0}^c}\big]&\le 2\E\big[\E_0[\|u(t_m)\|_H^2]\mathbf{1}_{\mathcal{G}_{N,c_0}^c}\big]+2\delta^{-2}\E[\mathbf{1}_{\mathcal{G}_{N,c_0}^c}].
\end{align*}
By \eqref{eq:sharp-u-moments}, \eqref{ineq:tail2}, and the definition of weak-$L^p$ norm tail bound given in \eqref{eq:sharp-budgets}, for $\varsigma<1/2$,
\begin{align*}
   \mathbb{S}_2
   &\le C(1+\big\|\|U_0\|_{\mathcal{H}_1^1}\big\|_{4,\infty}^4)N^{-1}+C\delta^{-2}N^{-2}\big\|\|U_0\|_{\mathcal{H}_1^1}\big\|_{4,\infty}^4\le C(1+\big\|\|U_0\|_{\mathcal{H}_1^1}\big\|_{4,\infty}^4)N^{-2\varsigma}.
\end{align*}

Combining the estimates for terms $\mathbb{S}_1$ and $\mathbb{S}_2$ yields the desired assertion for $t_m < \ell_N$:
\begin{align*}
&~~~\|u(t_m)-u^{N,\delta}_m\|_{L^2(\Omega;H)}^2\le  \mathbb{S}_1+\mathbb{S}_2\le C_{a,\varsigma}\big(1+\mathcal{J}_{\varsigma}+\big\|\|U_0\|_{\mathcal{H}_1^1}\big\|_{4,\infty}^4\big)\big(N^{-2\varsigma}+\tau(1+0\vee\log (1/\tau))\big).
\end{align*}

Finally, combining the estimates from Case 1, Case 2, and the initial time $m=0$, we conclude the proof of Theorem \ref{thm:sharp-main}.
\end{proof}

\begin{proof}[Proof of Remark~\ref{rem:sharp-Galerkin} (2)]
We first verify the initial approximation condition. Fix \(N\ge1\) and let \(M\to\infty\), so that \(\tau=T/M\to0\). Since \(u_0^{N,\delta}=\mathcal R_{N,N^{-\varsigma}}(\Pi_NU_0)\) is independent of \(\tau\), the assumed error estimate at \(m=0\) yields
\[
\mathbb E\bigl[\|u_0-u_0^{N,\delta}\|_H^2\bigr]
\le C_*N^{-2\varsigma}.
\]
Hence, it follows that
\[
\mathcal J_\varsigma
=
\sup_{N\ge1}
N^{2\varsigma}\,
\mathbb E\bigl[\|u_0-u_0^{N,\delta}\|_H^2\bigr]
\le C_*.
\]

We next derive the weak-\(L^4\) tail bound for \(\|u_0\|_H\). At time \(t_1=\tau\), the estimate $\|u_1^{N,\delta}\|_H\le\delta^{-1}\le N^\varsigma$ implies
\[
\bigl(0\vee(\|u(\tau)\|_H-N^\varsigma)\bigr)^2
\le
\|u(\tau)-u_1^{N,\delta}\|_H^2.
\]
Keeping \(N\) fixed and letting \(\tau\downarrow0\), the path continuity of \(u\) at \(0\), Fatou's lemma, and the assumed convergence estimate give
\begin{equation}\label{59}
\mathbb E\Bigl[\bigl(0\vee(\|u_0\|_H-N^\varsigma)\bigr)^2\Bigr]
\le C_*N^{-2\varsigma}.
\end{equation}
For sufficiently large \(R\), choose $N=\lfloor(R/2)^{1/\varsigma}\rfloor.$ Then, for some constant \(c_\varsigma>0\), we have $c_\varsigma R\le N^\varsigma\le R/2.$ Consequently, on the event \(\{\|u_0\|_H>R\}\),
\[
\bigl(0\vee(\|u_0\|_H-N^\varsigma)\bigr)^2\ge R^2/4.
\]
Using \eqref{59}, we obtain
\[
\mathbb P(\|u_0\|_H>R)
\le
\frac{4}{R^2}
\mathbb E\Bigl[\bigl(0\vee(\|u_0\|_H-N^\varsigma)\bigr)^2\Bigr]
\le
4C_*R^{-2}N^{-2\varsigma}
\le
C_\varsigma C_*R^{-4}.
\]
The bounded range of \(R\) is absorbed into the constant, and therefore
\[
\bigl\|\|u_0\|_H\bigr\|_{4,\infty}^4
\le C_\varsigma(1+C_*)<\infty.
\]
The proof is completed.
\end{proof}

\section{Numerical experiments}

We investigate the spatial and temporal convergence behavior of the proposed Cole--Hopf-based fully discrete scheme. The numerical experiments are conducted for \(\alpha=0\) and \(\alpha=1/8\), representing the space--time white-noise and rougher-noise regimes, respectively. All errors are measured in the root mean-square sense at the terminal time \(T=1\) against coupled numerical reference solutions.

We consider three representative initial conditions \(u_0\in H\), together with their corresponding Cole--Hopf transforms \(U_0\in\mathcal H_1^1\), satisfying the assumptions of Theorem \ref{thm:sharp-main}: 
\begin{align}\label{ineq:condition}
  J_\varsigma+\bigl\|\|U_0\|_{\mathcal H_1^1}\bigr\|_{4,\infty}^4<\infty.  
\end{align}
For each \(\alpha\in\{0,1/8\}\), the three initial data satisfy \eqref{ineq:condition} for every \(0<\varsigma<1/2-2\alpha\). In addition, they are chosen to represent different levels of the spatial regularity, ranging from a smooth constant profile to a rough initial datum, and thereby to examine the robustness of the proposed schemes with respect to the initial condition. Specifically, we consider the following initial data:
\[
\begin{aligned}
&\text{constant:}
&&u_0=0,\qquad U_0=1,\\
&\text{sine:}
&&u_0(x)=\sin(\pi x),\qquad
U_0(x)=\exp\left(\frac{1-\cos(\pi x)}{\pi}\right),\\
&\text{rough:}
&&u_0(x)=-\pi\sum_{j=1}^{\infty}\frac{\sin(j\pi x)}{j},\\
&&&U_0(x)
=
\exp\left(
\sum_{j=1}^{\infty}\frac{\cos(j\pi x)}{j^2}
-
\min_{y\in[0,1]}
\sum_{j=1}^{\infty}\frac{\cos(j\pi y)}{j^2}
\right).
\end{aligned}
\]
For the numerical implementation of the rough initial condition, the above series are truncated at \(j_{\max}=4096\). The resulting truncation error is negligible at the scale of experiments, since
$
\sum_{j>j_{\max}}j^{-2}=\mathcal{O}(j_{\max}^{-1}).
$  

We compute the root mean-square spatial and temporal errors for both the transformed approximation \(\{U_m^N\}_{m\in\{0,\ldots,M\}}\), defined by (16), and the reconstructed approximation \(\{u_m^{N,\delta}\}_{m\in\{0,\ldots,M\}}\), defined by \eqref{eq:sharp-output}. The following parameter settings are used to estimate the convergence rates:
\begin{itemize}
    \item Spatial convergence: $N_{\text{ref}} = 2^{10}$ and $\tau_{\text{ref}} = 2^{-11}$, with varying parameter $N = 2^i$ for $i = 3, \dots, 8$.
    \item Temporal convergence: $N_{\text{ref}} = 2^{9}$ and $\tau_{\text{ref}} = 2^{-12}$, with varying parameter $\tau = 2^{-i}$ for $i = 5, \dots, 10$.
\end{itemize}
All expectations are approximated by Monte Carlo averages over \(10^4\) independent realizations, and the same three initial conditions are used throughout.

\begin{figure}[!ht]
\centering
\includegraphics[height=5.0cm,width=10.5cm]{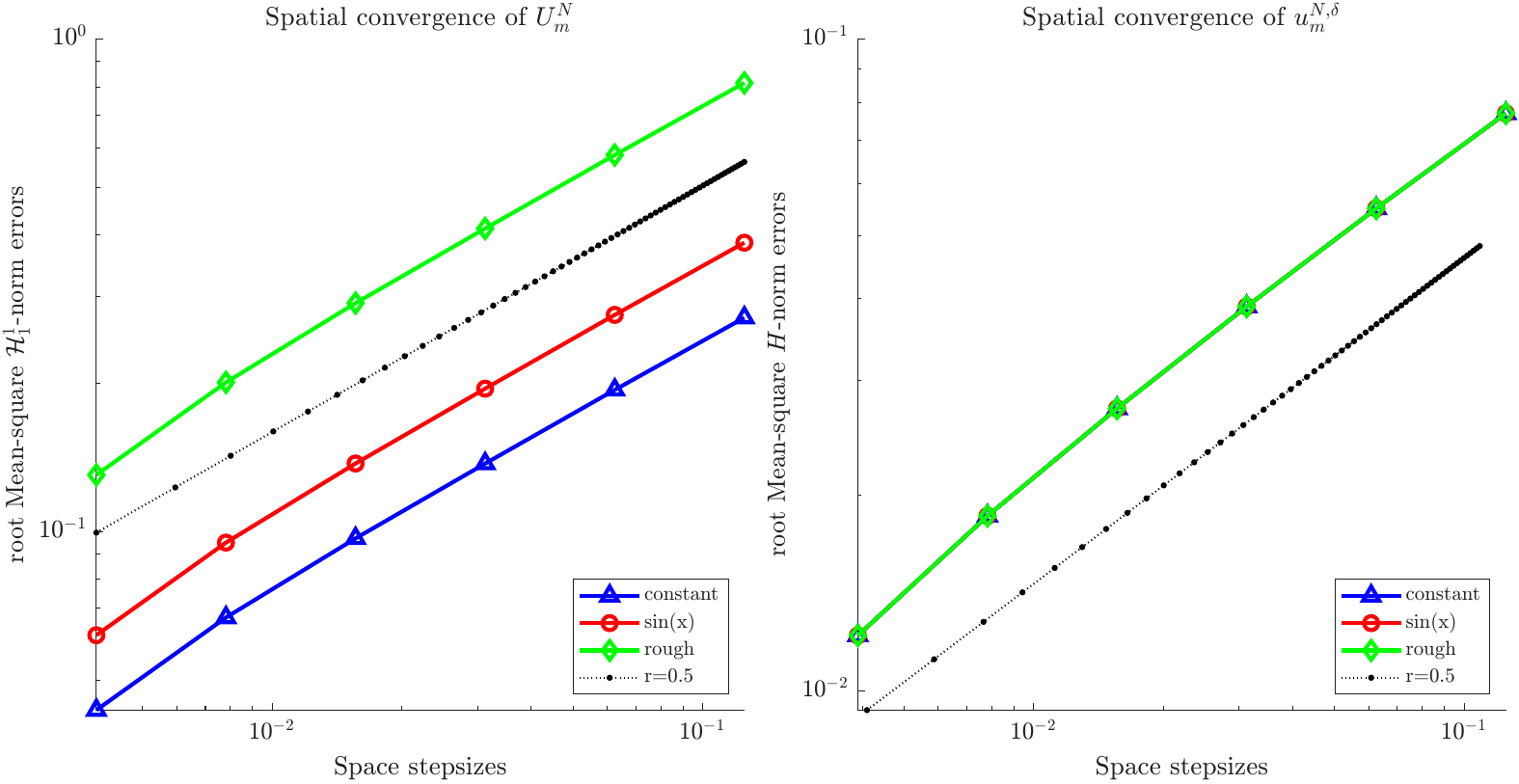}

\vspace{0.6cm}

\includegraphics[height=4.5cm,width=10.5cm]{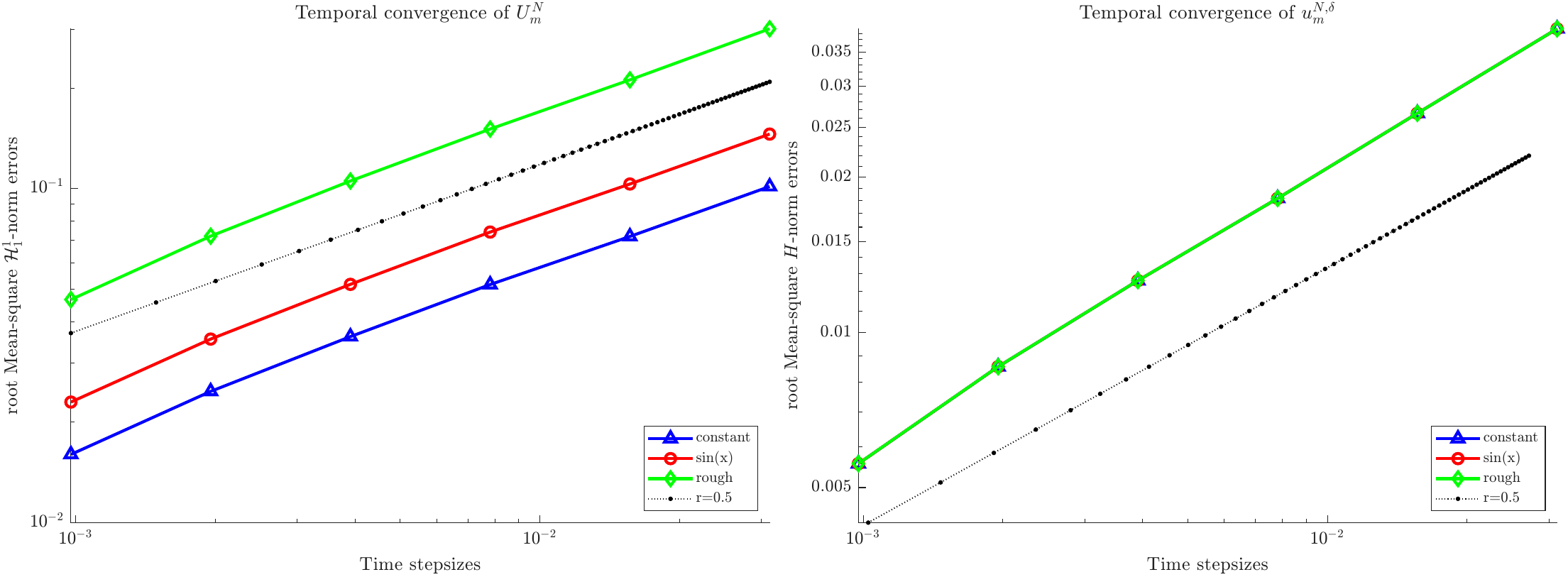}
\caption{
Spatial convergence results (top) and temporal convergence results (bottom)
for the transformed and reconstructed approximations in the case $\alpha=0$.
}\label{fig1}
\end{figure}

We first consider the space--time white noise case (\(\alpha=0\)). Figure \ref{fig1} reports the spatial and temporal root mean-square errors for both the transformed approximation $\{U^N_m\}_{m\in\{0,1,\ldots,M\}}$ and the reconstructed approximation $ \{u^{N,\delta}_m\}_{m\in\{0,1,\ldots,M\}}$. The results are displayed on log--log scales, together with reference lines of slope \(1/2\). For the transformed equation, the numerical errors exhibit approximately half-order decay with respect to both the spatial and temporal discretization parameters, consistently with the corresponding estimates \eqref{eq:sharp-linear-error} with $\gamma=1$ in Lemma~\ref{lem:sharp-linear-error}. 
The numerical approximation $ \{u^{N,\delta}_m\}_{m\in\{0,1,\ldots,M\}}$ shows essentially the same behavior: both the spatial and temporal errors decay with slopes close to \(1/2\) for all three initial data, as predicted by Theorem~\ref{thm:sharp-main}.



\begin{figure}[!ht]
\centering
\includegraphics[height=4.5cm,width=10.5cm]{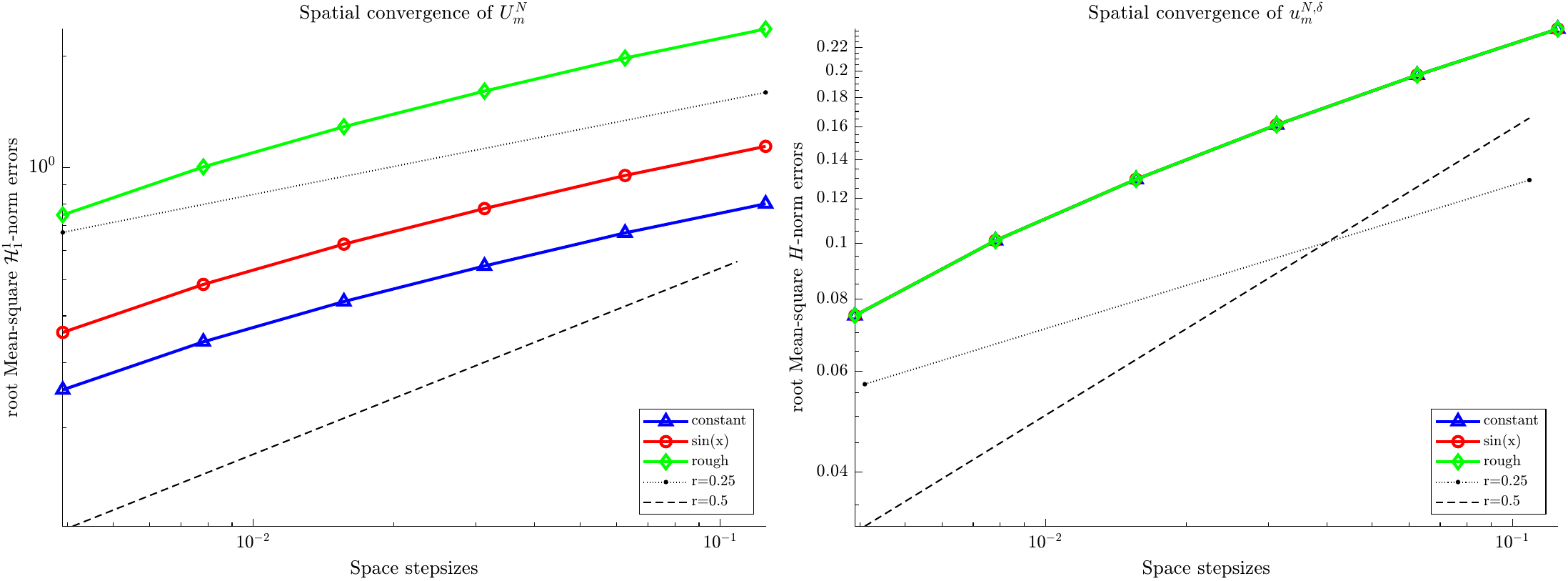}

\vspace{0.6cm}

\includegraphics[height=4.5cm,width=10.5cm]{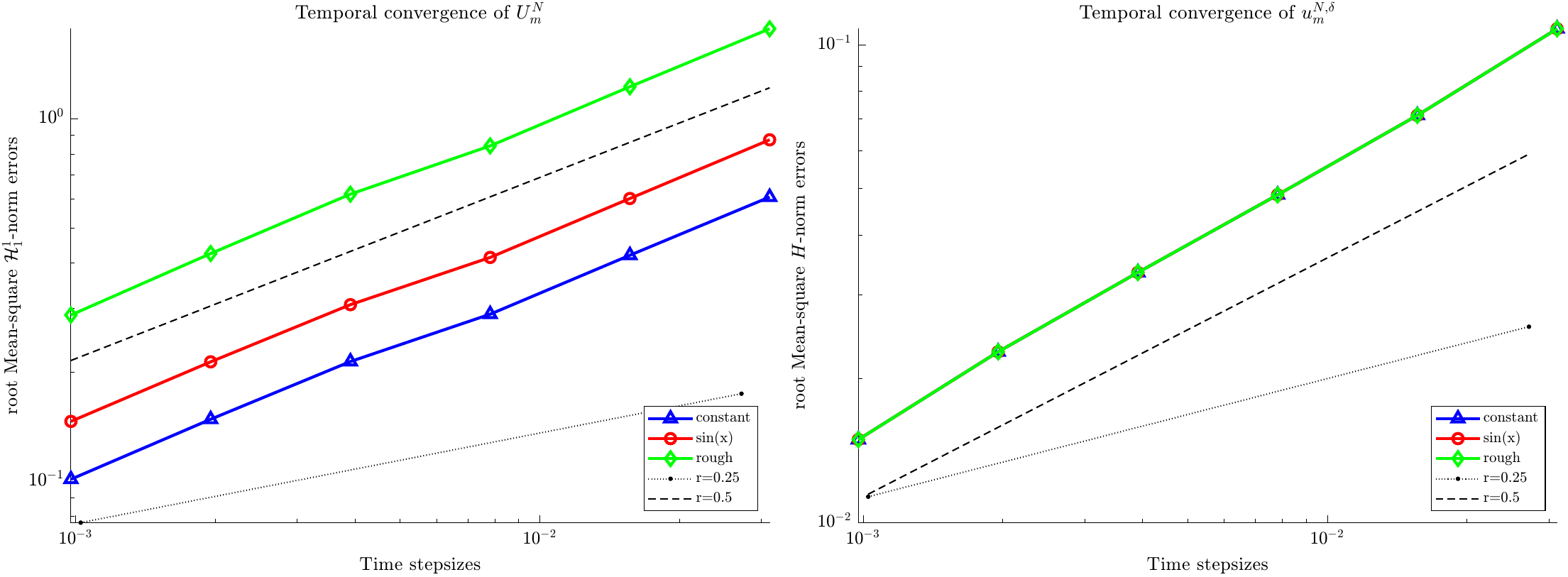}
\caption{
Spatial convergence results (top) and temporal convergence results (bottom)
for the transformed and reconstructed approximations in the case $\alpha=1/8$.
}\label{fig2}
\end{figure}



We next consider the rougher-noise case \(\alpha=1/8\), for which $1/2-2\alpha=1/4.$ According to Theorem \ref{thm:sharp-main}, the reconstructed approximation $ \{u^{N,\delta}_m\}_{m\in\{0,1,\ldots,M\}}$ admits spatial strong convergence of any order \(\varsigma<1/4\), whereas the temporal convergence rate remains of order \(1/2\), up to a logarithmic correction. The corresponding error bound reads
\[
\sup_{0\le m\le M}
\big(\E\big[\|u(t_m)-u_m^{N,\delta}\|_{H}^2\big]\big)^{1/2}
\le
C_{\varsigma,T}\bigl(N^{-\varsigma}+\tau^{1/2}(1+0\vee\ln(1/\tau))^{1/2}\bigr),\qquad 0<\varsigma<1/4.
\]

Figure \ref{fig2} also shows the spatial and temporal root mean-square errors for both the transformed approximation $\{U^N_m\}_{m\in\{0,1,\ldots,M\}}$ and the reconstructed approximation $ \{u^{N,\delta}_m\}_{m\in\{0,1,\ldots,M\}}$. In the spatial discretization, the observed slopes are close to \(1/4\) for all three initial data, in agreement with the regularity threshold \(1/2-2\alpha\). Compared with the space--time white-noise case (\(\alpha=0\)), where the spatial rate is close to \(1/2\), the rougher noise therefore leads to a clear deterioration of the spatial convergence rate.

In contrast, the temporal errors exhibit slopes close to \(1/2\) for both the transformed and reconstructed approximations. Thus, increasing the noise roughness from \(\alpha=0\) to \(\alpha=1/8\) substantially reduces the attainable spatial rate while leaving the observed temporal rate essentially unchanged. This behavior is consistent with estimates \eqref{eq:sharp-linear-error} with $\gamma=1$ in Lemma~\ref{lem:sharp-linear-error} and Theorem \ref{thm:sharp-main}, and numerically illustrates the distinct mechanisms governing the spatial and temporal discretization errors.

\appendix
\section{Analytic estimates}\label{app:analytic}

\begin{lemma}\label{lem:analytic}
For $0\le r\le\gamma$, $t\in(0,T]$, and $h>0$,
\begin{align*}
 \|S_1(t)\|_{\mathcal L(\HH_1^r;\HH_1^\gamma)}
       \le C_{\gamma,r,T}t^{-(\gamma-r)/2},\quad \|(S_1(h)-\mathrm I)S_1(t)\|_{\mathcal L(\HH_1^r)}
       \le C\min\{1,h/t\}.
\end{align*}
For $\gamma>1/2$, $\|f\|_\infty\le C_\gamma\|f\|_{\HH_1^\gamma}$.
Moreover,
\[
 \|(\mathrm I-\Pi_N)f\|_\infty
       \le CN^{-1/2}\|f\|_{\HH_1^1},\qquad
 \|f\|_{\HH_1^\gamma}
       \le\|f\|_H^{1-\gamma}\|f\|_{\HH_1^1}^\gamma
       \quad(0\le\gamma\le1).
\]
\end{lemma}
\begin{proof}
The semigroup estimates follow from standard fractional-power
estimates for analytic semigroups applied to $\Lambda_1=I-A_1$;
see \cite[Chapter 2, Theorem 6.13]{Pazy1983}.

Next, for $\gamma>1/2$, the Sobolev embedding gives $\|f\|_\infty\le C_\gamma\|f\|_{\HH_1^\gamma}$. Moreover,
\begin{align*}
 \|(\mathrm I-\Pi_N)f\|_\infty
 &\le\sqrt2\sum_{k>N}|\langle f,h_k\rangle_H|\le\sqrt2\Bigl(\sum_{k>N}(1+\lambda_k)^{-1}\Bigr)^{1/2}
       \|f\|_{\HH_1^1}\le  CN^{-1/2}\|f\|_{\HH_1^1}.
\end{align*}
For $0<\gamma<1$, H\"older's inequality yields
\begin{align*}
 \|f\|_{\HH_1^\gamma}^2
\le\Bigl(\sum_{k\ge0}|\langle f,h_k\rangle_H|^2\Bigr)^{1-\gamma}
       \Bigl(\sum_{k\ge0}(1+\lambda_k)|\langle f,h_k\rangle_H|^2\Bigr)^\gamma
 =\|f\|_H^{2(1-\gamma)}\|f\|_{\HH_1^1}^{2\gamma}.
\end{align*}
The cases $\gamma=0,1$ are immediate.
\end{proof}

\begin{lemma}\label{lem:sharp-multiplier}
For $0<r<1/2-2\alpha$, multiplication $M_{b_{k,\alpha}}v=b_{k,\alpha}v$ satisfies
\begin{equation}\label{eq:sharp-multiplier}
 \sum_{k\ge1}\|M_{b_{k,\alpha}}\|_{\mathcal L(\HH_1^r)}^2<\infty,
 \qquad
 \|\rho_\alpha v\|_{\HH_1^r}
 +\|G_\alpha(v)\|_{\mathcal L_2(\ell^2;\HH_1^r)}
       \le C_{\alpha,r}\|v\|_{\HH_1^r}.
\end{equation}
The bounds remain valid after applying $\Pi_N$, uniformly in $N$.
\end{lemma}
\begin{proof}
Let $\widehat v_j$ be the Fourier coefficients of $v$. Then
\[
 \|v\|_{\HH_1^r}^2
 \asymp\sum_{j\in\mathbb Z}(1+j^2)^r|\widehat v_j|^2,
 \qquad
 \widehat{\bigl(\cos(k\pi\cdot)v\bigr)}_j
 =\frac{\widehat v_{j-k}+\widehat v_{j+k}}2.
\]
Since $(1+|j\pm k|^2)^{r/2}
 \le C_r\bigl((1+j^2)^{r/2}+k^r\bigr),$ we have
\begin{align*}
 \|v\cos(k\pi\cdot)\|_{\HH_1^r}
 &\le C_r\Bigl[
 \Bigl(\sum_{j\in\mathbb Z}(1+j^2)^r|\widehat v_{j-k}|^2\Bigr)^{1/2}
 +\Bigl(\sum_{j\in\mathbb Z}(1+j^2)^r|\widehat v_{j+k}|^2\Bigr)^{1/2}
 \Bigr]\le C_r\bigl(\|v\|_{\HH_1^r}+k^r\|v\|_H\bigr).
\end{align*}
Recall that $b_{k,\alpha}=a_k(1-\cos(k\pi\cdot))$ and
$a_k=\sqrt2\,\pi^{2\alpha-1}k^{2\alpha-1}$. It holds that
\begin{align*}
 \|M_{b_{k,\alpha}}v\|_{\HH_1^r}
 &\le |a_k|\bigl(\|v\|_{\HH_1^r}
                 +\|v\cos(k\pi\cdot)\|_{\HH_1^r}\bigr)\le C_{\alpha,r}k^{2\alpha-1}(1+k^r)\|v\|_{\HH_1^r},
\end{align*}
and hence
\[
 \sum_{k\ge1}\|M_{b_{k,\alpha}}\|_{\mathcal L(\HH_1^r)}^2
 \le C_{\alpha,r}\sum_{k\ge1}k^{4\alpha-2}(1+k^{2r})
 <\infty,
 \qquad 4\alpha+2r<1.
\]
We obtain
\begin{align*}
 \|G_\alpha(v)\|_{\mathcal L_2(\ell^2;\HH_1^r)}^2\le\Bigl(\sum_{k\ge1}
       \|M_{b_{k,\alpha}}\|_{\mathcal L(\HH_1^r)}^2\Bigr)
       \|v\|_{\HH_1^r}^2,\quad
 \|\rho_\alpha v\|_{\HH_1^r}
 \le\frac12\Bigl(\sum_{k\ge1}
       \|M_{b_{k,\alpha}}\|_{\mathcal L(\HH_1^r)}^2\Bigr)
       \|v\|_{\HH_1^r}.
\end{align*}
Finally, $\|\Pi_N\|_{\mathcal L(\HH_1^r)}\le1$ implies the bounds are valid uniformly in $N$.
\end{proof}

\begin{lemma}\label{lem:volterra}
Let $\kappa\in(0,1)$. Suppose that \(F,g\in L^1_{\text{loc}}([0,T],\mathbb R_+)\). If 
$$F(t)\le g(t)+C\int_0^t(t-s)^{\kappa-1}F(s)\,ds,\qquad t\in[0,T],$$ 
then
\begin{equation}\label{cont vo}
 F(t)\le g(t)+C_{\kappa,C,T}
                 \int_0^t(t-s)^{\kappa-1}g(s)\,ds\qquad t\in[0,T].
\end{equation}
Moreover, let $M\in\mathbb N_+$, and set $\tau:=\frac{T}{M},~
t_m:=m\tau$, $m=0,1,\ldots,M$. 
 If $K\ge0$ and a nonnegative sequence $\{E_m\}_{m\in\{0,1,\ldots,M\}}$ satisfies $E_0\le K$ and, for \(m=1,\ldots,M\),
\begin{equation}\label{dis volt}
 E_m\le K+C\int_0^{t_m}(t_m-s)^{\kappa-1}E_{\lfloor s/\tau\rfloor }\,ds,
\end{equation}
then
$\max_{0\le m\le M}E_m\le C_{\kappa,C,T}K$ with $C_{\kappa,C,T}$ independent of $\tau$.
\end{lemma}

\begin{proof}
By the generalized weakly singular Gronwall inequality in
\cite{ye2007generalized}, we obtain
\[
F(t)\le g(t)+\int_0^t
\sum_{n=1}^{\infty}
\frac{(C\Gamma(\kappa))^{n}}{\Gamma(n\kappa)}
(t-s)^{n\kappa-1}g(s)\,ds .
\]
Moreover, for $0<r\le T$
\begin{align*}
\sum_{n=1}^{\infty}
\frac{(C\Gamma(\kappa))^{n}}{\Gamma(n\kappa)}
r^{n\kappa-1}=
r^{\kappa-1}\sum_{n=1}^{\infty}
\frac{(C\Gamma(\kappa))^{n}}{\Gamma(n\kappa)}
r^{(n-1)\kappa}  \le C_{\kappa,C,T}r^{\kappa-1} .
\end{align*}
Combining these estimates gives \eqref{cont vo}.

For the discrete inequality, choose $d>0$ such that
$Cd^\kappa/\kappa\le1/2$. Since
\[
 \int_0^{t_m}(t_m-s)^{\kappa-1}E_{\lfloor s/\tau\rfloor}\,ds
 =\sum_{j=0}^{m-1}E_j
 \int_{t_j}^{t_{j+1}}(t_m-s)^{\kappa-1}\,ds,
\]
for $\ell<m\le n$ with $t_n-t_\ell\le d$,
\begin{align*}
 E_m
 &\le K+C\int_0^{t_\ell}(t_m-s)^{\kappa-1}E_{\lfloor s/\tau\rfloor}\,ds
 +C\int_{t_\ell}^{t_m}(t_m-s)^{\kappa-1}E_{\lfloor s/\tau\rfloor}\,ds\\
 &\le K+\frac{CT^\kappa}{\kappa}\max_{0\le j\le\ell}E_j
 +\frac{C(t_m-t_\ell)^\kappa}{\kappa}
       \max_{\ell\le j<m}E_j\\
 &\le K+\frac{CT^\kappa}{\kappa}\max_{0\le j\le\ell}E_j
 +\frac12\max_{\ell\le j<m}E_j.
\end{align*}
Taking the maximum over $m$ yields 
$$\max_{\ell<m\le n}E_m\le K+\Bigl(\frac{CT^\kappa}{\kappa}+\frac12\Bigr)
\max_{0\le j\le\ell}E_j +\frac12\max_{\ell<m\le n}E_m.$$
 Notice that $\displaystyle\max _{0\le j \leq n} E_j=\max \{\max_{\ell<m\le n}E_m,\max_{0\le j\le\ell}E_j \}$, thus we have
 \begin{equation}\label{E11}
     \max_{0\le j\le n}E_j\le2K+\Bigl(\frac{2CT^\kappa}{\kappa}+1\Bigr)
       \max_{0\le j\le\ell}E_j.
\end{equation}
If $\tau\le d/2$, there are at most $2T/d+1$ blocks, thus using \eqref{E11} iteratively gives
\[
 \max_{0\le m\le M}E_m\le C_{\kappa,C,T}K.
\]
If $\tau>d/2$, then $M=T/\tau<2T/d$, and
\[
 E_m\le K+\frac{CT^\kappa}{\kappa}\max_{0\le j<m}E_j
 \quad\Longrightarrow\quad
 \max_{0\le j\le m}E_j
 \le K+\Big(1+\frac{CT^\kappa}{\kappa}\Big)\max_{0\le j<m}E_j.
\]
Since the number of time steps is bounded by $2T/d$, iterating this inequality yields $\max_{0\le m\le M}E_m\le C_{\kappa,C,T}K$, which completes the proof.
\end{proof}

\begin{lemma}\label{lem:heat-tail}
For $f\in L^1(0,1)$ and $0<t\le T$,
\begin{equation}\label{eq:heat-tail}
 \|S(t)\partial_xf\|_H\le Ct^{-3/4}\|f\|_1,
 \qquad
 \|(\mathrm I-P_N)S(t)\partial_xf\|_H
       \le Ct^{-3/4}e^{-cN^2t}\|f\|_1.
\end{equation}
For $v\in C([0,T];H)$, the Bochner integral
$\mathcal B(v)(t)=\int_0^tS(t-s)\partial_x(v(s)^2)\,ds$
belongs to $C([0,T];H)$, and
\begin{equation}\label{eq:B-local-Lipschitz}
 \|\mathcal B(v)-\mathcal B(w)\|_{C([0,T];H)}
 \le CT^{1/4}\|v-w\|_{C([0,T];H)}
        (\|v\|_{C([0,T];H)}+\|w\|_{C([0,T];H)}).
\end{equation}
\end{lemma}
\begin{proof}
Since $\|e_k'\|_\infty\le\sqrt2\pi k$ and $\sum_{k\ge1}k^2e^{-\pi^2k^2t} \le Ct^{-3/2},$ Parseval's identity yields
\begin{align*}
 \|S(t)\partial_xf\|_H^2
 &=\Big\|\sum_{k\ge1}e^{-\lambda_kt}\langle f,e_k'\rangle_{H}e_k\Big\|_H^2\le\|f\|_1^2\sum_{k\ge1}e^{-2\lambda_kt}\|e_k'\|_\infty^2
 \le2\pi^2\|f\|_1^2\sum_{k\ge1}k^2e^{-2\pi^2k^2t}
 \le Ct^{-3/2}\|f\|_1^2.
\end{align*}
Similarly, $e^{-2\pi^2k^2t}
\le e^{-\pi^2N^2t}e^{-\pi^2k^2t}$ for $k>N$ proves the second estimate in \eqref{eq:heat-tail}
\begin{align*}
 \|(\mathrm I-P_N)S(t)\partial_xf\|_H^2 \le2\pi^2e^{-\pi^2N^2t}\|f\|_1^2
       \sum_{k>N}k^2e^{-\pi^2k^2t}
 \le Ct^{-3/2}e^{-\pi^2N^2t}\|f\|_1^2.
\end{align*}
The Cauchy--Schwarz inequality gives
\begin{align}\label{eq:sharp-square-difference}
 \|v(s)^2-w(s)^2\|_1
 &\le\|v(s)-w(s)\|_H(\|v(s)\|_H+\|w(s)\|_H).
\end{align}
Hence $s\mapsto v(s)^2$ is continuous from $[0,T]$ to $L^1(0,1)$, and
\begin{align*}
 \int_0^t\|S(t-s)\partial_x(v(s)^2)\|_H\,ds
 &\le C\int_0^t(t-s)^{-3/4}\|v(s)\|_H^2\,ds\nonumber\le4Ct^{1/4}\|v\|_{C([0,T];H)}^2<\infty .
\end{align*}
Thus $\mathcal B(v)(t)$ is a well-defined $H$-valued Bochner integral and
\begin{align}\label{eq:sharp-B-size}
 \|\mathcal B(v)(t)\|_H
 \le4Ct^{1/4}\|v\|_{C([0,T];H)}^2 .
\end{align}
Moreover, \eqref{eq:heat-tail} and \eqref{eq:sharp-square-difference} imply \eqref{eq:B-local-Lipschitz} by the following
\begin{align*}
 &\|\mathcal B(v)(t)-\mathcal B(w)(t)\|_H\le C\int_0^t(t-s)^{-3/4}\|v(s)^2-w(s)^2\|_1\,ds\\
 \le& 4CT^{1/4}\|v-w\|_{C([0,T];H)}
       \bigl(\|v\|_{C([0,T];H)}+\|w\|_{C([0,T];H)}\bigr).
\end{align*}
For continuity, let $0<\varepsilon<t$ and $t-\varepsilon<r<t$.
Splitting the integrals at $t-\varepsilon$ gives
\begin{align*}
 &\|\mathcal B(v)(t)-\mathcal B(v)(r)\|_H\\
 &\le\bigl\|\bigl(S(\varepsilon)-S(\varepsilon-t+r)\bigr)
       \mathcal B(v)(t-\varepsilon)\bigr\|_H\\
 &\quad+C\int_{t-\varepsilon}^t(t-s)^{-3/4}\|v(s)\|_H^2\,ds
 +C\int_{t-\varepsilon}^r(r-s)^{-3/4}\|v(s)\|_H^2\,ds\\
 &\le\bigl\|\bigl(S(\varepsilon)-S(\varepsilon-t+r)\bigr)
       \mathcal B(v)(t-\varepsilon)\bigr\|_H
 +C\varepsilon^{1/4}\|v\|_{C([0,T];H)}^2 .
\end{align*}
The strong continuity of $S$ therefore yields
\[
 \limsup_{r\uparrow t}
 \|\mathcal B(v)(t)-\mathcal B(v)(r)\|_H
 \le C\varepsilon^{1/4}\|v\|_{C([0,T];H)}^2
 \longrightarrow0
 \qquad(\varepsilon\downarrow0).
\]
The same trick for $r\downarrow t$ proves the right continuity.
Finally, taking $t\to0$ in \eqref{eq:sharp-B-size} completes the proof.
\end{proof}

\begin{lemma}\label{lem:Burgers-rough-wellposed}
Let $0\le\alpha<\frac14$, and let $u_0$ be an
$\mathcal F_0$-measurable $H$-valued random variable. Then
\eqref{eq:sharp-sbe} admits a unique adapted mild solution with paths in
$C([0,T];H)$
\begin{equation}\label{mild u}
    u(t)=S(t) u_0+\int_0^t S(t-s) \partial_x\left(u(s)^2\right) d s+\int_0^t S(t-s) (-A)^\alpha d W(s),\qquad a.s.
\end{equation}
The uniqueness is understood up to
indistinguishability.
\end{lemma}
\begin{proof}
The existence is supplied by the construction in the proof of
Lemma~\ref{lem:sharp-existence}  and  \eqref{eq:sharp-mild}.
That construction first solves the linear equation
\eqref{eq:sharp-pam}, establishes its strict positivity,
and then applies Lemma~\ref{lem:finite-noise-log-chain}
to obtain $u=\partial_xU/U\in C([0,T];H)$.

For completeness, write
$Z_\alpha(t):=\int_0^tS(t-s)(-A)^\alpha\,dW(s)$
for the stochastic convolution in \eqref{mild u}.
Then
\[
 \mathbb E[\|Z_\alpha(t)\|_H^2]
 =\frac12\sum_{k\ge1}\lambda_k^{2\alpha-1}
       \bigl(1-e^{-2\lambda_kt}\bigr)
 \le C\sum_{k\ge1}k^{4\alpha-2}<\infty .
\]

For $0\le s<t\le T$, independence of the increments gives
\begin{align*}
 &\mathbb E[\|Z_\alpha(t)-Z_\alpha(s)\|_H^2]=\frac12\sum_{k\ge1}\lambda_k^{2\alpha-1}
 \Bigl[1-e^{-2\lambda_k(t-s)}
 +(1-e^{-\lambda_k(t-s)})^2(1-e^{-2\lambda_ks})\Bigr].
\end{align*}
If $0<\eta<1/2-2\alpha$, then
$1-e^{-2x}\le Cx^\eta$ and $(1-e^{-x})^2\le Cx^\eta$ for $x\ge0,$ and hence
\begin{align*}
 \mathbb E[\|Z_\alpha(t)-Z_\alpha(s)\|_H^2]
 &\le C(t-s)^\eta
       \sum_{k\ge1}\lambda_k^{2\alpha-1+\eta}\le C(t-s)^\eta\sum_{k\ge1}k^{4\alpha-2+2\eta}
 \le C(t-s)^\eta .
\end{align*}
Since $Z_\alpha(t)-Z_\alpha(s)$ is Gaussian,
\[
 \mathbb E[\|Z_\alpha(t)-Z_\alpha(s)\|_H^p]
 \le C_p\bigl(\mathbb E[\|Z_\alpha(t)-Z_\alpha(s)\|_H^2]\bigr)^{p/2}
 \le C_p|t-s|^{p\eta/2}.
\]
Choosing $p>2/\eta$, Kolmogorov's continuity theorem gives a
continuous $H$-valued modification of $Z_\alpha$.

For every $u\in C([0,T];H)$, Lemma~\ref{lem:heat-tail} yields
\begin{align*}
 \int_0^t\|S(t-s)\partial_x(u(s)^2)\|_H\,ds
 &\le C\int_0^t(t-s)^{-3/4}\|u(s)\|_H^2\,ds\le4Ct^{1/4}\|u\|_{C([0,T];H)}^2<\infty ,
\end{align*}
so the nonlinear convolution is a well-defined, continuous
$H$-valued Bochner integral. Together with the strong continuity of
$S$, this gives a mild solution  in $C([0,T];H)$.

Finally, if $u$ and $v$ are two such solutions with the same initial
value and noise, then
\begin{align*}
 &\sup_{r\le t}\|u(r)-v(r)\|_H\le C\bigl(\|u\|_{C([0,T];H)}+\|v\|_{C([0,T];H)}\bigr)
 \int_0^t(t-\tau)^{-3/4}
 \sup_{r\le\tau}\|u(r)-v(r)\|_H\,d\tau .
\end{align*}
Lemma~\ref{lem:volterra} with $\kappa=\frac14$ proves the pathwise uniqueness.
\end{proof}

\begin{lemma}\label{lem:deterministic-HC}
Let $V_0,W_0\in\HH_1^1$ have the strictly positive minima and set
\[
 v(t)=\frac{\partial_x(S_1(t)V_0)}{S_1(t)V_0},\qquad
 w(t)=\frac{\partial_x(S_1(t)W_0)}{S_1(t)W_0}.
\]
Then $v,w\in C([0,T];H)$ solve the deterministic Burgers equation,
their $H$ norms do not increase, and
\begin{equation}\label{eq:deterministic-comparison}
 \|v(t)-P_Nw(t)\|_H
 \le\|v(0)-P_Nw(0)\|_H
       +Ct^{1/4}(\|v(0)\|_H^2+\|w(0)\|_H^2).
\end{equation}
\end{lemma}
\begin{proof}
The maximum principle gives $S_1(t)V_0\ge\min_{[0,1]}V_0>0$
and $S_1(t)W_0\ge\min_{[0,1]}W_0>0$. Moreover, {the strong continuity of \(S_1\) on \(\mathcal H_1^1\), together with the boundedness of
\(\partial_x\colon\mathcal H_1^1\to H\), gives}
\[
 S_1(\cdot)V_0,\ S_1(\cdot)W_0\in C([0,T];\HH_1^1),
 \qquad
 \partial_x(S_1(\cdot)V_0),\ \partial_x(S_1(\cdot)W_0)\in C([0,T];H),
\]
and hence $v,w\in C([0,T];H)$.

For $t>0$, $\partial_tS_1(t)V_0=\partial_{xx}S_1(t)V_0$ and
$\partial_xS_1(t)V_0|_{x=0,1}=0$. The direct differentiation gives
\[
 v_t
 =\partial_x\left(\frac{\partial_{xx}(S_1(t)V_0)}{S_1(t)V_0}\right)
 =\partial_{xx}v+\partial_{x}(v^2),
 \qquad v(t,0)=v(t,1)=0,
\]
and similarly for $w$. Consequently,
\begin{align*}
 \frac12\frac{d}{dt}\|v(t)\|_H^2
 &=\langle v(t),\partial_{xx}v(t)\rangle_H
       +\langle v(t),\partial_{x}(v(t)^2)\rangle_H=-\|\partial_{x}v(t)\|_H^2 .
\end{align*}
Thus, for $0<\varepsilon<t$,
\[
 \frac12\bigl(\|v(t)\|_H^2-\|v(\varepsilon)\|_H^2\bigr)
 =-\int_\varepsilon^t\|\partial_{x}v(s)\|_H^2\,ds\le0,
\]
and letting $\varepsilon\downarrow0$ proves
$\|v(t)\|_H\le\|v(0)\|_H$. The same calculation gives
$\|w(t)\|_H\le\|w(0)\|_H$.

The facts that $S(t)P_N=P_NS(t)$ and
$\|P_N\|_{\mathcal L(H)}\le1$ yield
\begin{align*}
 &\|v(t)-P_Nw(t)\|_H\\
 &\le\|S(t)(v(0)-P_Nw(0))\|_H
       +\|\mathcal B(v)(t)\|_H+\|P_N\mathcal B(w)(t)\|_H\\
 &\le\|v(0)-P_Nw(0)\|_H
 +C\int_0^t(t-s)^{-3/4}
       \bigl(\|v(s)\|_H^2+\|w(s)\|_H^2\bigr)\,ds\\
 &\le\|v(0)-P_Nw(0)\|_H
 +4Ct^{1/4}\bigl(\|v(0)\|_H^2+\|w(0)\|_H^2\bigr),
\end{align*}
where Lemma~\ref{lem:heat-tail} was used for the two convolutions.
\end{proof}

\section{The Cole--Hopf identities}\label{app:colehopf}

\begin{lemma}[Finite-noise Feynman--Kac formula]
\label{lem:FeynmanKac}
For $K\ge1$, let $\hat{U}^K$ solve \eqref{eq:sharp-pam} with $b_{k,\alpha}$ set
to zero for $k>K$ and drift
$\rho_{\alpha,K}=\frac12\sum_{k\le K}b_{k,\alpha}^2$.
Let $R^x$ be reflected Brownian motion on $[0,1]$ with generator
$\partial_{xx}$ (cf. \cite{lions1984stochastic}), independent of $\mathcal F_T$, and let $\E_R$
denote expectation over $R^x$. Then
\begin{equation}\label{eq:feynman}
 \hat{U}^K(t,x)=\E_R\bigg[U_0(R^x(t))
    \exp\bigg\{\sum_{k\le K}\int_0^tb_{k,\alpha}(R^x(t-s))\,d\beta_k(s)
        \bigg\}\bigg].
\end{equation}
For each fixed $(t,x)$ the identity holds almost surely. The
right-hand side has a jointly continuous modification for which
the equality holds simultaneously for all $(t,x)$.
\end{lemma}

\begin{proof}
The following computation follows
\cite[Sect.~5]{bertini1994stochastic} with respect to $R^x(t)$. Its transition density $p_t$ satisfies
\begin{align}\label{Neumann-semigroup}
 \E_R[\varphi(R^x(t))]
 =\int_0^1p_t(x,y)\varphi(y)\,dy
 =(S_1(t)\varphi)(x), \quad \varphi \in L^\infty(0,1).
\end{align}
Conditional on $\mathcal F_0$ and the trajectories of $R^x$, the stochastic
integral in \eqref{eq:feynman} is Gaussian, and
\begin{align*}
 &\E\bigg[\exp\bigg\{a\sum_{k\le K}\int_0^t
 b_{k,\alpha}(R^x(t-s))\,d\beta_k(s)\bigg\}
 \bigg|\,\mathcal F_0,R^x\bigg]\\
 &\qquad=\exp\bigg\{\frac{a^2}{2}\sum_{k\le K}\int_0^t
 b_{k,\alpha}(R^x(t-s))^2\,ds\bigg\}
 \le\exp\bigg\{\frac{a^2t}{2}\sum_{k\le K}\|b_{k,\alpha}\|_\infty^2\bigg\}.
\end{align*}
Let $V_K(t,x)$ denote the right-hand side of \eqref{eq:feynman}.
Jensen's inequality gives
\begin{align}\label{eq:finite-noise-FK-moment}
 \E_0[|V_K(t,x)|^2]
 &\le\|U_0\|_\infty^2
 \exp\bigg\{2t\sum_{k\le K}\|b_{k,\alpha}\|_\infty^2\bigg\}\le C_{K,T}\|U_0\|_{\HH_1^1}^2 .
\end{align}
Fix $t,x$ and define
\[
 M_s^{t,x}
 =\sum_{k\le K}\int_0^s b_{k,\alpha}(R^x(t-z))\,d\beta_k(z),
 \qquad 0\le s\le t .
\]
It\^o's formula gives
\begin{align}\label{eq:finite-noise-exponential}
 de^{M_s^{t,x}}
 =e^{M_s^{t,x}}\bigg[
 \sum_{k\le K}b_{k,\alpha}(R^x(t-s))\,d\beta_k(s)
 +\rho_{\alpha,K}(R^x(t-s))\,ds\bigg].
\end{align}
For bounded Borel $\varphi$, the Markov property at time $t-s$ yields
\begin{align}\label{eq:finite-noise-Markov}
 &\E_R\left[U_0(R^x(t))e^{M_s^{t,x}}
       \varphi(R^x(t-s))\right]\nonumber\\
 &=\int_0^1p_{t-s}(x,y)\varphi(y)
 \E_R\bigg[U_0(R^y(s))
 \exp\bigg\{\sum_{k\le K}\int_0^s
 b_{k,\alpha}(R^y(s-r))\,d\beta_k(r)\bigg\}\bigg]dy\nonumber\\
 &=\int_0^1p_{t-s}(x,y)\varphi(y)V_K(s,y)\,dy
 =\bigl(S_1(t-s)[\varphi V_K(s)]\bigr)(x).
\end{align}
Integrating \eqref{eq:finite-noise-exponential}, multiplying by
$U_0(R^x(t))$, taking $\E_R$, and using
\eqref{eq:finite-noise-Markov} give
\begin{align*}
 V_K(t,x)
 &=(S_1(t)U_0)(x)
 +\int_0^t\bigl(S_1(t-s)[\rho_{\alpha,K}V_K(s)]\bigr)(x)\,ds+\sum_{k\le K}\int_0^t
 \bigl(S_1(t-s)[b_{k,\alpha}V_K(s)]\bigr)(x)\,d\beta_k(s).
\end{align*}
Thus $V_K$ and $\hat{U}^K$ satisfy the same linear mild equation in $H$.
Semigroup contractivity and It\^o's isometry imply
\begin{align*}
 &\E_0\big[\|V_K(t)-\hat{U}^K(t)\|_H^2\big]\le\Big(2T\|\rho_{\alpha,K}\|_\infty^2
 +2\sum_{k\le K}\|b_{k,\alpha}\|_\infty^2\Big)
 \int_0^t\E_0\big[\|V_K(s)-\hat{U}^K(s)\|_H^2\big]\,ds,
\end{align*}
and hence Gronwall's inequality gives
$\E_0\big[\|V_K(t)-\hat{U}^K(t)\|_H^2\big]=0$.

For the pointwise identification, the Neumann heat kernel satisfies
\[
 \sup_{x\in[0,1]}\|p_\tau(x,\cdot)\|_H^2
 =\sup_{x\in[0,1]}\sum_{j\ge0}e^{-2\lambda_j\tau}|h_j(x)|^2
 \le2\sum_{j\ge0}e^{-2\pi^2j^2\tau}
 \le C_T\tau^{-1/2}.
\]
Subtracting the pointwise mild equations and using Cauchy--Schwarz
and It\^o's isometry therefore give
\begin{align*}
 &\E_0\big[|V_K(t,x)-\hat{U}^K(t,x)|^2\big]\\
 &\le\Big(2T\|\rho_{\alpha,K}\|_\infty^2
 +2\sum_{k\le K}\|b_{k,\alpha}\|_\infty^2\Big)
 \int_0^t\sup_{y\in[0,1]}\|p_{t-s}(y,\cdot)\|_H^2
 \E_0\big[\|V_K(s)-\hat{U}^K(s)\|_H^2\big]\,ds\\
 &\le C_{K,T}\int_0^t(t-s)^{-1/2}
 \E_0\big[\|V_K(s)-\hat{U}^K(s)\|_H^2\big]\,ds=0.
\end{align*}
The identity at $t=0$ is immediate. Finally,
\[
 \sup_{x\in[0,1]}|\hat{U}^K(t,x)-\hat{U}^K(s,x)|
 \le C\|\hat{U}^K(t)-\hat{U}^K(s)\|_{\HH_1^1},
\]
so the continuous $\HH_1^1$-valued version of $\hat{U}^K$ is jointly
continuous in $(t,x)$ and is the asserted continuous modification
of $V_K$.
\end{proof}

\begin{lemma}\label{lem:finite-noise-log-chain}
Let $U\in C([0,T],\HH_1^1)$ be a mild solution of \eqref{eq:sharp-pam} and $\inf_{(t,x)\in[0,T]\times[0,1]}U(t,x)>0$ a.s. Then $g=\log U\in C([0,T],\HH_1^1)$ {a.s.,} and satisfies \eqref{eq:sharp-log-weak} simultaneously for every $\varphi\in\HH_1^1$ and $t\in[0,T]$.
\end{lemma}
\begin{proof}
Since $\partial_xg=\partial_xU/U$,
\begin{equation*}
 \|g(t)-g(s)\|_{\HH_1^1}\le C\left(\frac1{\displaystyle\inf_{(t,x)\in[0,T]\times[0,1]}U(t,x)}
 +\frac{\displaystyle\sup_{0\le t\le T}\|U(t)\|_{\HH_1^1}}
        {\big(\displaystyle\inf_{(t,x)\in[0,T]\times[0,1]}U(t,x)\big)^2}\right)
 \|U(t)-U(s)\|_{\HH_1^1},
\end{equation*}
so $g\in C([0,T];\HH_1^1)$ a.s. The mild equation also gives
\[
 dU=(A_1U+\rho_\alpha U)\,dt
       +\sum_kUb_{k,\alpha}\,d\beta_k
 \qquad\text{in }(\HH_1^1)^*,
\]
where
\begin{align*}
 {|{}_{(\mathcal H_1^1)^*}\!\left\langle A_1v,\psi\right\rangle_{\mathcal H_1^1}|}
 &=|\langle \partial_xv,\partial_x\psi\rangle_H|
 \le\|v\|_{\HH_1^1}\|\psi\|_{\HH_1^1},\\
 \|A_1v+\rho_\alpha v\|_{(\HH_1^1)^*}
 &\le C\|v\|_{\HH_1^1},\qquad
 \sum_k\|vb_{k,\alpha}\|_H^2
 \le\Bigl(\sum_k\|b_{k,\alpha}\|_\infty^2\Bigr)\|v\|_H^2
 \le C_\alpha\|v\|_H^2.
\end{align*}

For fixed \(\varphi\in\mathcal H_1^1\) and \(R>1\), choose $\log_R \in C_b^2(\mathbb{R})$ such that
$\log_R(z)=\ln z$ for $z\in[R^{-1},R]$, and define $F_{R,\varphi}\colon H\to \mathbb R$ by 
\[
F_{R,\varphi}(v):=\int_0^1 \varphi(x)\,\log_R(v(x))\,dx,
\qquad v\in H.
\]
Then $DF_{R,\varphi}(v)[h]=\int_0^1 \varphi(x)\,(\log_R)'(v(x))h(x)\,dx$,  $D^2F_{R,\varphi}(v)[h,h]=
\int_0^1 \varphi(x)\,(\log_R)''(v(x))h(x)^2\,dx,
$ and
\begin{align*}
|DF_{R,\varphi}(v)[h]|
\le
C_{R,\varphi}\|h\|_H,\qquad
|D^2F_{R,\varphi}(v)[h,h]|
\le
C_{R,\varphi}\|h\|_H^2
\end{align*}
{with both derivatives jointly continuous in $(v,h)\in H\times H$.}
Moreover, for $v\in\HH^1_1$,
\[
 \nabla F_{R,\varphi}(v)=\varphi(\log_R)'(v)\in\HH_1^1,\qquad
 \|\nabla F_{R,\varphi}(v)\|_{\HH_1^1}
 \le C_{R,\varphi}(1+\|v\|_{\HH_1^1}),
\]
and $\nabla F_{R,\varphi}:\HH_1^1\to\HH_1^1$ is continuous. These properties follow from the embedding $\mathcal{H}_1^1 \hookrightarrow C([0,1])$ and the boundedness of $(\log_R)^{\prime}$ and $(\log_R)^{\prime \prime}$. 

Define
\[
 \tau_R:=\inf\left\{t\in[0,T]:
 \min_{x\in[0,1]}U(t,x)\le R^{-1}\ \text{or}\
 \max_{x\in[0,1]}U(t,x)\ge R\ \text{or}\
 \|U(t)\|_{\HH_1^1}\ge R\right\}\wedge T .
\]
It holds that
\[
 \frac{\mathbf1_{\{s<\tau_R\}}}{R^2\max\{1,C_\alpha\}}
 \sum_k\|U(s)b_{k,\alpha}\|_H^2
 \le\frac{C_\alpha}{\max\{1,C_\alpha\}}\le1 .
\]
Thus the variational It\^o formula
\cite[Theorem~3.1]{krylov2013ito} applies to
\[
 \frac{U(\cdot\wedge\tau_R)}{R\max\{1,\sqrt{C_\alpha}\}},
 \qquad
 v\longmapsto F_{R,\varphi}\bigl(R\max\{1,\sqrt{C_\alpha}\}\,v\bigr).
\]
On the event $\Big\{\min_{x\in[0,1]}U_0(x)>R^{-1},\
 \|U_0\|_\infty<R,\ \|U_0\|_{\HH_1^1}<R\Big\},$ the stopped trajectory lies in $[R^{-1},R]$, where $\log_R=\ln$.
Up to $\tau_R$, It\^o's formula therefore yields
\begin{align*}
 d\langle g,\varphi\rangle_H
 &=\bigg[-\int_0^1\partial_xU\,\partial_x\!\left(\frac{\varphi}{U}\right)dx
 +\int_0^1\rho_\alpha\varphi\,dx
 -\frac12\sum_k\int_0^1b_{k,\alpha}^2\varphi\,dx\bigg]dt+\sum_k\langle b_{k,\alpha},\varphi\rangle_H\,d\beta_k
\end{align*}
Here notice that $\rho_\alpha=\frac12\sum_kb_{k,\alpha}^2$ and
\begin{align*}
 -\int_0^1\partial_xU\,\partial_x\!\left(\frac{\varphi}{U}\right)dx
 &=-\int_0^1\frac{\partial_xU}U\partial_x\varphi\,dx
 +\int_0^1\left|\frac{\partial_xU}U\right|^2\varphi\,dx=-\langle\partial_xg,\partial_x\varphi\rangle_H
 +\int_0^1|\partial_xg|^2\varphi\,dx.
\end{align*}
Consequently, on the same event,
\begin{align*}
 \langle g(t\wedge\tau_R),\varphi\rangle_H
 &=\langle\log U_0,\varphi\rangle_H+\int_0^{t\wedge\tau_R}
 \Big[-\langle\partial_xg,\partial_x\varphi\rangle_H
 +\int_0^1|\partial_xg|^2\varphi\,dx\Big]ds+\sum_k\langle b_{k,\alpha},\varphi\rangle_H
       \beta_k(t\wedge\tau_R).
\end{align*}

Almost surely, $\inf_{(t,x)\in[0,T]\times[0,1]}U(t,x)>0, \sup_{(t,x)\in[0,T]\times[0,1]}U(t,x)<\infty$ and $\sup_{0\le t\le T}\|U(t)\|_{\HH_1^1}<\infty,$ so the event holds and $\tau_R=T$ for all sufficiently large
integers $R$. Letting $R\to\infty$ proves the integrated identity for
every $t\in[0,T]$ and the fixed $\varphi$.

Finally,
\[
 \mathbb E\Big[\sup_{t\le T}
 \Big|\sum_{k>K}\langle b_{k,\alpha},\varphi\rangle_H\beta_k(t)\Big|^2\Big]
 \le4T\|\varphi\|_H^2\sum_{k>K}\|b_{k,\alpha}\|_H^2
 \longrightarrow0\qquad(K\to\infty),
\]
which gives a continuous version of the series. All other terms
are continuous in $(t,\varphi)\in[0,T]\times\HH_1^1$. A countable
dense subset of $\HH_1^1$ therefore makes the identity simultaneous
in $t$ and $\varphi$.
\end{proof}

\section*{Acknowledgement}
The authors acknowledge the use of AI-assisted tools for writing
assistance and for exploring calculations and experiment codes. All
mathematical content, verification, and final conclusions remain
the sole responsibility of the authors.

\end{document}